%
\newif\ifloadreferences\loadreferencestrue
%
%
%
%
%
\let\myfrac=\frac%
\input eplain %
\let\frac=\myfrac%
\let\myfootnote=\footnote%
\input amstex \input epsf %
\let\footnote=\myfootnote%
%
%
\loadeufm\loadmsam\loadmsbm\message{symbol names}\UseAMSsymbols\message{,}%
\magnification 1200 %
\font\myfontdefault=cmr10%
\newif\ifmakebiblio%
\newif\ifinappendices%
\newif\ifundefinedreferences%
\newif\ifchangedreferences%
\makebibliofalse%
\undefinedreferencesfalse%
\changedreferencesfalse%
%
%
%
%
%
\def\setcatcodes{\catcode`\!=0 \catcode`\\=11}%
{\global\let\noe=\noexpand%
\catcode`\@=11 \catcode`\_=11 \setcatcodes%
!global!def!_@@internal@@makeref#1{%
!global!expandafter!def!csname #1ref!endcsname##1{%
!csname _@#1@##1!endcsname%
!expandafter!ifx!csname _@#1@##1!endcsname!relax%
    !write16{#1 ##1 not defined - run saving references}%
    !undefinedreferencestrue%
!fi}}%
!global!def!_@@internal@@makelabel#1{%
!global!expandafter!def!csname #1label!endcsname##1{%
!edef!temptoken{!csname #1info!endcsname}%
!ifloadreferences%
    !expandafter!ifx!csname _@#1@##1!endcsname!relax%
        !write16{#1 ##1 not hitherto defined - rerun saving references}%
        !changedreferencestrue%
    !else%
        !expandafter!ifx!csname _@#1@##1!endcsname!temptoken%
        !else%
            !write16{#1 ##1 reference has changed - rerun saving references}%
            !changedreferencestrue%
        !fi%
    !fi%
!else%
    !expandafter!edef!csname _@#1@##1!endcsname{!temptoken}%
    !edef!textoutput{!write!references{\global\def\_@#1@##1{!temptoken}}}%
    !textoutput%
!fi}}%
!global!def!makecounter#1{!_@@internal@@makelabel{#1}!_@@internal@@makeref{#1}}%
!unsetcatcodes%
}
%
%
%
%
%
\def\turnintolatin#1{\ifcase #1 _\or i\or ii\or iii\or iv\or v\or vi\or vii\or viii\or ix\or x\or xi\or xii\or xiii\or xiv\or xv\or xvi\or xvii\or xviii\or xix\or xx\or xxi\or xxii\or xxiii\or xxiv\or xxv\or xxvi\fi}%
\def\alphanum#1{\ifcase #1 _\or A\or B\or C\or D\or E\or F\or G\or H\or I\or J\or K\or L\or M\or N\or O\or P\or Q\or R\or S\or T\or U\or V\or W\or X\or Y\or Z\fi}%
\newwrite\references%
\ifloadreferences{\catcode`\@=11 \catcode`\_=11 \global\def\_@citation@AgmonNirenberg{1}
\global\def\_@citation@Alexander{2}
\global\def\_@citation@AlmgrenSimon{3}
\global\def\_@citation@Andrews{4}
\global\def\_@citation@AndrewsLangfordMcCoy{5}
\global\def\_@citation@GrueterJost{6}
\global\def\_@citation@Hamilton{7}
\global\def\_@citation@HamiltonII{8}
\global\def\_@citation@Mantegazza{9}
\global\def\_@citation@MacSmiMT{10}
\global\def\_@citation@MaximoNunesSmith{11}
\global\def\_@citation@PacardXu{12}
\global\def\_@citation@RobbinSalamon{13}
\global\def\_@citation@RosSmiDT{14}
\global\def\_@citation@Rudin{15}
\global\def\_@citation@Schneider{16}
\global\def\_@citation@Schwarz{17}
\global\def\_@citation@SmiEC{18}
\global\def\_@citation@SmiEMCFI{19}
\global\def\_@citation@SmiEMCFII{20}
\global\def\_@citation@Tromba{21}
\global\def\_@citation@Weber{22}
\global\def\_@citation@WhiteI{23}
\global\def\_@citation@WhiteII{24}
\global\def\_@citation@Ye{25}
\global\def\_@head@Introduction{1}
\global\def\_@subhead@PrescribedCurvatureProblems{1.1}
\global\def\_@eqn@AreaMinusVolumeFunctional{\relax \unhbox \voidb@x \hbox {{\relax \tenrm (1)}}}
\global\def\_@eqn@PrescribedMeanCurvature{\relax \unhbox \voidb@x \hbox {{\relax \tenrm (2)}}}
\global\def\_@proc@MainExistenceResult{1.1.1}
\global\def\_@subhead@MorseHomology{1.2}
\global\def\_@eqn@GradientFlowOperator{\relax \unhbox \voidb@x \hbox {{\relax \tenrm (3)}}}
\global\def\_@eqn@LinearisedGradientFlowOperator{\relax \unhbox \voidb@x \hbox {{\relax \tenrm (4)}}}
\global\def\_@subhead@EternalForcedMeanCurvatureFlows{1.3}
\global\def\_@eqn@ForcedMCFEquation{\relax \unhbox \voidb@x \hbox {{\relax \tenrm (5)}}}
\global\def\_@eqn@ControllingFunction{\relax \unhbox \voidb@x \hbox {{\relax \tenrm (6)}}}
\global\def\_@eqn@Subcritical{\relax \unhbox \voidb@x \hbox {{\relax \tenrm (7)}}}
\global\def\_@eqn@LittleLambdaAndLambda{\relax \unhbox \voidb@x \hbox {{\relax \tenrm (8)}}}
\global\def\_@eqn@AdmissableElements{\relax \unhbox \voidb@x \hbox {{\relax \tenrm (9)}}}
\global\def\_@proc@Separation{1.3.1}
\global\def\_@proc@Compactness{1.3.2}
\global\def\_@proc@EndPoints{1.3.3}
\global\def\_@proc@ConcentrationEllipticCase{1.3.4}
\global\def\_@proc@ConcentrationParabolicCase{1.3.5}
\global\def\_@subhead@ConstructingTheMorseHomology{1.4}
\global\def\_@subhead@Acknowledgements{1.5}
\global\def\_@head@CompactnessAndEndPoints{2}
\global\def\_@subhead@GeometricBounds{2.1}
\global\def\_@proc@DoesNotContainBigSpheres{2.1.1}
\global\def\_@proc@IsNotContainedInSmallSpheres{2.1.2}
\global\def\_@proc@CConvexity{2.1.3}
\global\def\_@proc@BoundedDiameter{2.1.4}
\global\def\_@subhead@Compactness{2.2}
\global\def\_@proc@QuasiMaximumLemma{2.2.1}
\global\def\_@proc@BoundedCurvature{2.2.2}
\global\def\_@subhead@EndPoints{2.3}
\global\def\_@proc@UniqueContinuation{2.3.1}
\global\def\_@proc@EitherConstantOrOpen{2.3.2}
\global\def\_@proc@EndPointsExist{2.3.3}
\global\def\_@head@Admissability{3}
\global\def\_@subhead@Pinching{3.1}
\global\def\_@proc@GeneralisedSimonsIdentity{3.1.1}
\global\def\_@proc@DerivativeOfPrincipalCurvatures{3.1.2}
\global\def\_@proc@PreservationOfPinching{3.1.3}
\global\def\_@proc@Hopf{3.1.4}
\global\def\_@subhead@NonCollapsing{3.2}
\global\def\_@proc@PreservationOfNonCollapsing{3.2.1}
\global\def\_@head@SingularPerturbation{4}
\global\def\_@subhead@AsymptoticSeries{4.1}
\global\def\_@eqn@TaylorSeriesOfH{\relax \unhbox \voidb@x \hbox {{\relax \tenrm (10)}}}
\global\def\_@eqn@TaylorSeriesOfN{\relax \unhbox \voidb@x \hbox {{\relax \tenrm (11)}}}
\global\def\_@eqn@FlowFunction{\relax \unhbox \voidb@x \hbox {{\relax \tenrm (12)}}}
\global\def\_@eqn@DefinitionOfPsi{\relax \unhbox \voidb@x \hbox {{\relax \tenrm (13)}}}
\global\def\_@eqn@AsymptoticSeries{\relax \unhbox \voidb@x \hbox {{\relax \tenrm (14)}}}
\global\def\_@eqn@AsymptoticSeriesOfPsi{\relax \unhbox \voidb@x \hbox {{\relax \tenrm (15)}}}
\global\def\_@subhead@SphericalHarmonicsAndFormalSolutions{4.2}
\global\def\_@eqn@DefinitionOfL{\relax \unhbox \voidb@x \hbox {{\relax \tenrm (16)}}}
\global\def\_@eqn@DefinitionOfP{\relax \unhbox \voidb@x \hbox {{\relax \tenrm (17)}}}
\global\def\_@proc@FormalSolutions{4.2.1}
\global\def\_@eqn@InductiveDefinitionOfPhi{\relax \unhbox \voidb@x \hbox {{\relax \tenrm (18)}}}
\global\def\_@eqn@InductiveDefinitionOfX{\relax \unhbox \voidb@x \hbox {{\relax \tenrm (19)}}}
\global\def\_@subhead@TheKappaDependentRescaledTimeVariable{4.3}
\global\def\_@proc@UniformNormOfK{4.3.2}
\global\def\_@eqn@BoundsForA{\relax \unhbox \voidb@x \hbox {{\relax \tenrm (20)}}}
\global\def\_@eqn@BoundsForH{\relax \unhbox \voidb@x \hbox {{\relax \tenrm (21)}}}
\global\def\_@proc@UniformBoundsOnCoefficients{4.3.3}
\global\def\_@subhead@ExactSolutions{4.4}
\global\def\_@proc@BoundedInverse{4.4.1}
\global\def\_@proc@SingularPerturbations{4.4.2}
\global\def\_@head@Concentration{5}
\global\def\_@subhead@AnotherRescaledTimeVariable{5.1}
\global\def\_@eqn@DefinitionOfTildePhi{\relax \unhbox \voidb@x \hbox {{\relax \tenrm (22)}}}
\global\def\_@subhead@Bootstrapping{5.2}
\global\def\_@eqn@DefinitionOfTildePsi{\relax \unhbox \voidb@x \hbox {{\relax \tenrm (23)}}}
\global\def\_@proc@Bootstrapping{5.2.1}
\global\def\_@subhead@RecoveringTheFlow{5.3}
\global\def\_@proc@RecoveringTheFlow{5.3.1}
\global\def\_@subhead@AsymptoticSeriesAgain{5.4}
\global\def\_@proc@FlowsAreAsymptoticToFormalSolution{5.4.1}
\global\def\_@head@WeaklySmoothManifolds{A}
\global\def\_@head@Bibliography{B}
 }%
\else{\openout\references=references.tex }%
\fi%
%
%
\newcount\headno%
\global\headno=0%
\def\headinfo{\ifinappendices\alphanum\headno\else\the\headno\fi}%
\def\nextheadno[#1]{\global\advance\headno by 1 \global\subheadno=0 \global\procno=0 \headinfo\headlabel{#1}}%
\makecounter{head}%
%
%
\newcount\subheadno%
\global\subheadno=0%
\def\subheadinfo{\headinfo.\the\subheadno}%
\def\nextsubheadno[#1]{\global\advance\subheadno by 1 \global\procno=0 \subheadinfo\subheadlabel{#1}}%
\makecounter{subhead}%
%
%
\newcount\procno%
\global\procno=0%
\def\procinfo{\subheadinfo.\the\procno}%
\def\nextprocno{\global\advance\procno by 1 \procinfo}%
\makecounter{proc}%
%
%
\newcount\figno%
\global\figno=0%
\def\figinfo{\subheadinfo.\the\figno}%
\def\nextfigno{\global\advance\figno by 1 \figinfo}%
\makecounter{fig}%
%
%
\newcount\eqnno%
\global\eqnno=0%
\def\eqninfo{\text{{\rm (\the\eqnno)}}}%
\def\nexteqnno[#1]{\global\advance\eqnno by 1\hbox to 0em{\eqnlabel{#1}}\eqninfo}%
\makecounter{eqn}%
%
%
%
%
%
\def\gobbleeight#1#2#3#4#5#6#7#8{}%
\newcount\citationno%
\global\citationno=0%
\def\citationinfo{\the\citationno}%
\makecounter{citation}%
\newwrite\biblio%
\def\newref#1#2{%
\def\temptext{#2}%
\edef\bibliotextoutput{\expandafter\gobbleeight\meaning\temptext}%
\global\advance\citationno by 1\citationlabel{#1}%
\ifmakebiblio%
    \edef\fileoutput{\write\biblio{\noindent\hbox to 0pt{\hss$[\the\citationno]$}\hskip 0.2em\bibliotextoutput\medskip}}%
    \fileoutput%
\fi}%
\def\cite#1{%
$[\citationref{#1}]$%
\ifmakebiblio%
    \edef\fileoutput{\write\biblio{#1}}%
    \fileoutput%
\fi%
}%
%
%
%
%
\let\mypar=\par%
\edef\Pagetitle={Blank}\headline={\hfil\Pagetitle\hfil}%
\edef\Pagefooter={Blank}\footline={\hfil\Pagefooter\hfil}%
%
%
\newcount\showpagenumflag%
\global\showpagenumflag=0 %
\def\nextoddpage%
{\newpage\ifodd\pageno%
\else\global\showpagenumflag=0 %
\null\vfil\eject%
\global\showpagenumflag=1 %
\fi}%
%
%
\font\headfont=cmb12%
\def\newhead#1[#2]%
{\ifhmode\mypar\fi%
\ifnum\headno=0 \else\goodbreak\bigskip\fi%
{\headfont\noindent\nextheadno[#2]\ - #1.}
\nobreak\medskip}%
%
%
\def\newsubhead#1[#2]%
{\ifhmode\mypar\fi%
\ifnum\subheadno=0 \else\goodbreak\medskip\fi%
{\bf\noindent\nextsubheadno[#2]\ - #1.\ }}%
%
%
\newif\ifinproclaim%
\global\inproclaimfalse%
\def\proclaim#1{%
\goodbreak\medskip
\bgroup\inproclaimtrue%
\noindent{\bf #1}%
\nobreak\medskip\sl}%
\def\noskipproclaim#1{%
\goodbreak\medskip%
\bgroup\inproclaimtrue%
\noindent{\bf #1}\nobreak\sl}%
\def\endproclaim{\mypar\egroup\nobreak\medskip\ignorespaces}%
%
%
%
\newcount\xpos\newcount\ypos
\def\makelabelgrid{%
\xpos=-5 \ypos=-5 %
\loop\ifnum\xpos<6 %
{\loop\ifnum\ypos<6 %
\def\labeltext{x}%
\ifnum\xpos=0\def\labeltext{+}\fi%
\ifnum\ypos=0\def\labeltext{+}\fi%
\placelabel[\xpos][\ypos]{\labeltext}%
\advance\ypos by 1 %
\repeat}%
\advance\xpos by 1 %
\repeat}%
\def\placelabel[#1][#2]#3{{%
\setbox10=\hbox{\raise #2cm \hbox{\hskip #1cm #3}}%
\ht10=0pt \dp10=0pt \wd10=0pt \box10}}%
%
%
%
%
\def\myitem#1{\noindent\hbox to .5cm{\hfill#1\hss}}%
%
%
%
%
%
%
%
%
%
\font\sansseriften=cmss10%
\font\sansserifseven=cmss7%
\font\sansseriffive=cmss5%
\newfam\sansseriffam%
\textfont\sansseriffam=\sansseriften%
\scriptfont\sansseriffam=\sansserifseven%
\scriptscriptfont\sansseriffam=\sansseriffive%
\def\mathsf{\fam\sansseriffam}%
%
%
%
\font\boldten=cmb10%
\font\boldseven=cmb7%
\font\boldfive=cmb5%
\newfam\mathboldfam%
\textfont\mathboldfam=\boldten%
\scriptfont\mathboldfam=\boldseven%
\scriptscriptfont\mathboldfam=\boldfive%
\def\mathbf{\fam\mathboldfam}%
%
%
%
\font\mycmmiten=cmmi10%
\font\mycmmiseven=cmmi7%
\font\mycmmifive=cmmi5%
\newfam\mycmmifam%
\textfont\mycmmifam=\mycmmiten%
\scriptfont\mycmmifam=\mycmmiseven%
\scriptscriptfont\mycmmifam=\mycmmifive%
\def\hexa#1{\ifcase #1 0\or 1\or 2\or 3\or 4\or 5\or 6\or 7\or 8\or 9\or A\or B\or C\or D\or E\or F\fi}%
\mathchardef\mathi="7\hexa\mycmmifam7B%
\mathchardef\mathj="7\hexa\mycmmifam7C%
%
%
\font\mymsbmten=msbm10 at 8pt%
\font\mymsbmseven=msbm7 at 5.6pt
\font\mymsbmfive=msbm5 at 4pt%
\newfam\mymsbmfam%
\textfont\mymsbmfam=\mymsbmten%
\scriptfont\mymsbmfam=\mymsbmseven%
\scriptscriptfont\mymsbmfam=\mymsbmfive%
\mathchardef\mybeth="7\hexa\mymsbmfam69%
\mathchardef\mygimmel="7\hexa\mymsbmfam6A%
\mathchardef\mydaleth="7\hexa\mymsbmfam6B%
%
%
%
%
\def\proof{{\noindent\bf Proof:\ }}%
\def\remark{{\noindent\bf Remark:\ }}%
\def\qed{~$\square$}%
\def\makeop#1{\global\expandafter\def\csname op#1\endcsname{{\text{{\rm #1}}}}}%
\def\makeopsmall#1{\global\expandafter\def\csname op#1\endcsname{{\text{{\rm \lowercase{#1}}}}}}%
%
%
\def\munion{\mathop{\cup}}%
\def\minter{\mathop{\cap}}%
%
%
\makeop{Ext}%
\makeop{Int}%
\makeop{Dist}%
\makeop{Diam}%
\makeop{Length}%
%
%
%
%
%
\def\mlim{\mathop{{\text{Lim}}}}%
\def\mlimsup{\mathop{{\text{LimSup}}}}%
\def\msup{\mathop{{\text{Sup}}}}%
\def\minf{\mathop{{\text{Inf}}}}%
%
%
\makeop{Dim}%
\makeop{Ker}%
\makeop{Coker}%
\makeop{Tr}%
\makeop{Adj}%
\makeop{Det}%
\makeop{End}%
\makeop{Lin}%
\makeop{Symm}%
\makeop{Mult}%
%
%
\makeop{dx}%
\makeop{dy}%
\makeop{dz}%
\makeop{dt}%
\makeop{dVol}%
\makeop{dArea}%
\makeop{Supp}%
\makeop{Hess}%
\makeop{Lip}%
%
%
\makeop{Re}%
\makeop{Im}%
\makeop{Arg}%
\makeop{Log}%
\makeop{Exp}%
%
%
\makeopsmall{Cos}%
\makeopsmall{Sin}%
\makeopsmall{Tan}%
\makeopsmall{Sec}%
\makeopsmall{Cosec}%
\makeopsmall{Cot}%
\makeopsmall{ArcCos}%
\makeopsmall{ArcSin}%
\makeopsmall{ArcTan}%
\makeopsmall{ArcSec}%
\makeopsmall{ArcCosec}%
\makeopsmall{ArcCot}%
%
%
\makeopsmall{Cosh}%
\makeopsmall{Sinh}%
\makeopsmall{Tanh}%
\makeopsmall{ArcCosh}%
\makeopsmall{ArcSinh}%
\makeopsmall{ArcTanh}%
%
%
\makeop{Vol}%
\makeop{Area}%
\makeop{Riem}%
\makeop{Ric}%
\makeop{Scal}%
\makeop{Euc}%
\makeop{Imm}%
\makeop{Emb}%
%
%
\makeop{Id}%
\makeop{Ad}%
\makeop{O}%
\makeop{SO}%
\makeop{SL}%
\makeop{GL}%
\makeop{Conf}%
\makeop{Homeo}%
\makeop{Diff}%
\makeop{Isom}%
%
%
\makeop{Ind}%
\makeop{Sig}%
\makeop{Spec}%
%
%
\makeop{Conv}%
\makeop{Max}%
\makeop{Min}%
\makeop{Mod}%
\makeop{Deg}%
\makeop{loc}%
%
%
%
%
%
%
%
%
%
%
%
%
%
 %
%
%
%
%
%
\makeop{Sup}%
\makeopsmall{Log}%
\makeop{Index}%
\makeop{bdd}%
\makeop{in}
\catcode`\@=11
\def\triplealign#1{\null\,\vcenter{\openup1\jot \m@th %
\ialign{\strut\hfil$\displaystyle{##}\quad$&\hfil$\displaystyle{{}##}$&$\displaystyle{{}##}$\hfil\crcr#1\crcr}}\,}
\def\multiline#1{\null\,\vcenter{\openup1\jot \m@th %
\ialign{\strut$\displaystyle{##}$\hfil&$\displaystyle{{}##}$\hfil\crcr#1\crcr}}\,}
\catcode`\@=12
\newref{AgmonNirenberg}{Agmon S., Nirenberg L., Lower bounds and uniqueness theorems for solutions of differential equations in a Hilbert space, {\sl Comm. Pure Appl. Math.}, {\bf 20}, (1967), 207–-229}
\newref{Alexander}{Alexander S., Locally convex hypersurfaces of negatively curved spaces, {\sl Proc. Amer. Math. Soc.}, {\bf 64}, (1977), no. 2, 321--325}
\newref{AlmgrenSimon}{Almgren F. J., Simon L., Existence of embedded solutions of Plateau's problem, {\sl Ann. Scuola Norm. Sup. Pisa Cl. Sci.}, (\bf 6), (1979), no. 3, 447--495}
\newref{Andrews}{Andrews B., Non-collapsing in mean-convex mean curvature flow, {\sl Geometry \& Topology}, {\bf 16}, (2012), 1413--1418}
\newref{AndrewsLangfordMcCoy}{Andrews B., Langford M., McCoy J., Non-collapsing in fully nonlinear curvature flows, {\sl Annales de l'Institut Henri Poincar\'e - Analyse non lin\'eaire}, {\bf 30}, (2013), no. 1, 23--32}
\newref{GrueterJost}{Gr\"uter M., Jost J., On embedded minimal disks in convex bodies. {\sl Ann. Inst. H. Poincar\'e, Anal. Non Lin\'eaire}, {\bf 3}, (1986), no. 5, 345--390}
\newref{Hamilton}{Hamilton R. S., Convex hypersurfaces with pinched second fundamental form, {\sl Comm. Anal. Geom.}, {\bf 2}, (1994), no. 1, 167--172}
\newref{HamiltonII}{Hamilton R. S., The inverse function theorem of Nash and Moser, {\sl Bull. Amer. Math. Soc.}, {\bf 7}, (1982), no. 1, 65--222}
\newref{Mantegazza}{Mantegazza C., {\sl Lecture notes on mean curvature flow}, Progress in Mathematics, {\bf 290}, Springer Verlag, Basel, (2011)}
\newref{MacSmiMT}{Macarini L., Smith G., Morse homology of forced mean curvature flows, {\sl in preparation}}
\newref{MaximoNunesSmith}{M\'aximo D., Nunes I., Smith G., Free boundary minimal annuli in convex three-manifolds, to appear in {\sl J. Diff. Geom}}
\newref{PacardXu}{Pacard F., Xu X., Constant mean curvature spheres in Riemannian manifolds, {\sl Manuscripta Math.}, {\bf 128}, (2009), no. 3, 275--295}
\newref{RobbinSalamon}{Robbin J., Salamon D., The spectral flow and the Maslov index, {\sl Bull. London Math. Soc.}, {\bf 27}, (1995), no. 1, 1--33}
\newref{RosSmiDT}{Rosenberg H., Smith G., Degree theory of immersed hypersurfaces, arXiv:1010.1879}
\newref{Rudin}{Rudin W., {\sl Principles of mathematical analysis}, International Series in Pure \& Applied Mathematics, McGraw-Hill, (1976)}
\newref{Schneider}{Schneider M., Closed magnetic geodesics on closed hyperbolic Riemann surfaces, {\sl Proc. London Math. Soc.}, (2012), {\bf 105}, 424--446}
\newref{Schwarz}{Schwarz M., {\sl Morse homology}, Progress in Mathematics, {\bf 111}, Birkh\"auser Verlag, Basel, (1993)}
\newref{SmiEC}{Smith G., Constant curvature hypersurfaces and the Euler characteristic, arXiv:1103.3235}
\newref{SmiEMCFI}{Smith G., Eternal forced mean curvature flows I - a compactness result, {\sl Geom. Dedicata}, {\bf 176}, no. 1, (2014), 11--29}
\newref{SmiEMCFII}{Smith G., Eternal Forced Mean Curvature Flows II - Existence, arXiv:1508.05688}
\newref{Tromba}{Tromba A, {\sl A theory of branched minimal surfaces}, Springer Monographs in Mathematics, Springer-Verlag, Berlin, Heidelberg, (2012)}
\newref{Weber}{Weber J., Morse homology for the heat flow, {\sl Math. Z.}, {\bf 275}, (2013), no. 1, 1--54}
\newref{WhiteI}{White B., Every three-sphere of positive Ricci curvature contains a minimal embedded torus, {\sl Bull. Amer. Math. Soc.}, {\bf 21}, (1989), no. 1, 71--75}
\newref{WhiteII}{White B., The space of minimal submanifolds for varying Riemannian metrics, {\sl Indiana Univ. Math. J.}, {\bf 40}, (1991), no. 1, 161--200}
\newref{Ye}{Ye R., Foliation by constant mean curvature spheres, {\sl Pacific J. Math.}, {\bf 147}, (1991), no. 2, 381--396}
\catcode`\@=11
\def\triplealign#1{\null\,\vcenter{\openup1\jot \m@th %
\ialign{\strut\hfil$\displaystyle{##}$&$\displaystyle{{}##}\hfil$&$\displaystyle{{}##}$\hfil\crcr#1\crcr}}\,}
\def\multiline#1{\null\,\vcenter{\openup1\jot \m@th %
\ialign{\strut$\displaystyle{##}$\hfil&$\displaystyle{{}##}$\hfil\crcr#1\crcr}}\,}
\catcode`\@=12
\def\Pagetitle{\hfil}
\def\Pagefooter{\hfil{\myfontdefault\folio}\hfil}
\null\vfill
\def\centre{\rightskip=0pt plus 1fil \leftskip=0pt plus 1fil \spaceskip=.3333em \xspaceskip=.5em \parfillskip=0em \parindent=0em}%
\def\textmonth#1{\ifcase#1\or January\or Febuary\or March\or April\or May\or June\or July\or August\or September\or October\or November\or December\fi}
\font\abstracttitlefont=cmr10 at 14pt {\abstracttitlefont\centre
Eternal forced mean curvature flows III - Morse homology.\par}
\bigskip
{\centre 12th\ January 2016\par}
\bigskip
{\centre Graham Smith\par}
\bigskip
{\centre Instituto de Matem\'atica,\par
UFRJ, Av. Athos da Silveira Ramos 149,\par
Centro de Tecnologia - Bloco C,\par
Cidade Universit\'aria - Ilha do Fund\~ao,\par
Caixa Postal 68530, 21941-909,\par
Rio de Janeiro, RJ - BRASIL\par}
\bigskip
\noindent{\bf Abstract:}\ We complete the theoretical framework required for the construction of a Morse homology theory for certain types of forced mean curvature flows. The main result of this paper describes the asymptotic behaviour of these flows as the forcing term tends to infinity in a certain manner. This result allows the Morse homology to be explicitely calculated, and will permit us to show in forthcoming work that, for a large family of smooth positive functions, $F$, defined over a $(d+1)$-dimensional flat torus, there exist at least $2^{d+1}$ distinct, locally strictly convex, Alexandrov-embedded hyperspheres of mean curvature prescribed at every point by $F$.
%
%
\par
\vfill
\nextoddpage
\global\pageno=1
\myfontdefault
\def\Pagetitle{\hfil{\myfontdefault Eternal forced mean curvature flows III - Morse homology}\hfil}
\def\Pagefooter{\hfil{\myfontdefault\folio}\hfil}
\newhead{Introduction}[Introduction]
\newsubhead{Prescribed curvature problems}[PrescribedCurvatureProblems]
The problem of constructing hypersurfaces of constant curvature subject to geometric or topological restrictions is a standard one of riemannian geometry. The related problem of constructing hypersurfaces whose curvature is {\sl prescribed} by some function of the ambient space is essentially complementary, in that a full understanding of the one generally entails a full understanding of the other. However, as they are usually more straightforward, prescribed curvature problems often serve as a better testing ground for the development of new ideas, as will be the case here.
\par
Consider, therefore, the problem of constructing immersed hyperspheres of prescribed mean curvature inside a $(d+1)$-dimensional riemannian manifold, $M$.\footnote*{The convention will be adopted throughout this text that the mean curvature of an immersed hypersurface is equal to the {\sl arithmetic mean} of its principal curvatures, as opposed to their sum.} It is useful at this stage to introduce some terminology. Thus, $S^d$ will denote the unit sphere inside $\Bbb{R}^{d+1}$. A smooth immersion, $e:S^d\rightarrow M$, will be said to be an {\bf Alexandrov-embedding} whenever it extends to a smooth immersion, $\tilde{e}:\overline{B}^{d+1}\rightarrow M$, where $\overline{B}^{d+1}$ here denotes the closed unit ball inside $\Bbb{R}^{d+1}$. $\hat{\Cal{E}}(M)$ will denote the set of all smooth Alexandrov-embeddings from $S^d$ into $M$, and $\Cal{E}(M)$ will denote its quotient under the action by reparametrisation of the group of orientation-preserving, smooth diffeomorphisms of $S^d$. The spaces $\hat{\Cal{E}}(M)$ and $\Cal{E}(M)$ will be furnished respectively with the topology of $C^k$-convergence for all $k$, and its induced quotient topology, making $\Cal{E}(M)$, in particular, into a weakly smooth manifold (c.f. Appendix \headref{WeaklySmoothManifolds}). Finally, an element, $[e]$, of $\Cal{E}(M)$ will be referred to simply as an {\bf Alexandrov-embedding}, and will be identified, at times with a representative element, $e$, in $\hat{\Cal{E}}(M)$, and at times with its image in $M$.
\par
For the purposes of this paper, we may suppose that the extension of any Alexandrov-embedding is actually unique up to diffeomorphism.\footnote{${}^\dagger$}{Although this would seem to be the case whenever the ambient manifold is not homeomorphic to a sphere, we know of no such a result in the literature. In the present paper, we will only be concerned with locally strictly convex Alexandrov-embeddings in flat tori. Since all such maps lift to embeddings bounding convex sets in Euclidean space (c.f. \cite{Alexander}), uniqueness of the extension readily follows in this case.} Thus, for a smooth function, $F:M\rightarrow\Bbb{R}$, the functional, $\Cal{F}:\Cal{E}(M)\rightarrow\Bbb{R}$, will be defined by
$$
\Cal{F}([e]) := \opVol([e]) - d\int_{\overline{B}^{d+1}}(F\circ\tilde{e})\opdVol_{\tilde{e}},\eqnum{\nexteqnno[AreaMinusVolumeFunctional]}
$$
where $\tilde{e}:\overline{B}^{d+1}\rightarrow M$ here denotes the extension of $e$, and $\opdVol_{\tilde{e}}$ denotes the volume form that it induces over the closed, unit ball. The integral on the right-hand side will be viewed as a weighted volume of the interior of $e$, so that the functional, $\Cal{F}$, will be loosely referred to as the ``{\bf area-minus-volume}'' functional. The critical points of this functional are precisely those Alexandrov embedded hyperspheres whose mean curvatures are {\bf prescribed} at every point by the function, $F$; that is, those elements, $[e]$, of $\Cal{E}(M)$ such that
$$
H_e - F\circ e = 0,\eqnum{\nexteqnno[PrescribedMeanCurvature]}
$$
where $H_e$ here denotes the mean curvature of the Alexandrov-embedding, $e$.
\par
The fact that they arise as critical points suggests that Alexandrov-embedded hyperspheres of prescribed mean curvature should be amenable to study by differential-topological techniques, which should then yield information about their number, at least for generic data. This idea, which is far from new, has already been used by numerous authors to obtain some quite inspiring results in this, and related, settings (c.f., for example \cite{AlmgrenSimon}, \cite{GrueterJost}, \cite{Schneider}, \cite{Tromba}, \cite{WhiteI} and \cite{WhiteII}). Of particular relevance to the current discussion, however, is our own, relatively straightforward, result, \cite{RosSmiDT}, which shows that, under suitable circumstances, the number of critical points of $\Cal{F}$ is bounded below by the Euler characteristic of the ambient space. Although this yields existence in certain cases, we have found it rather unsatisfactory, as, on the one hand, it provides little new information in the case where the Euler characteristic of the ambient space vanishes - for example when the ambient space is $3$-dimensional - and, on the other, even when this Euler characteristic does not vanish, the number of solutions it yields generally falls fall short of what Morse theory would lead us to expect.
\par
Although Morse homology theory would appear to be the natural approach for obtaining the best possible lower bounds on the number of critical points of the functional, $\Cal{F}$, its development in the current setting has proven to be far from trivial, largely due to the technical challenges involved in making any progress in the theory of mean curvature flows. However, the remarkable non-collapsing theorems recently obtained by Ben Andrews et al. (c.f., for example, \cite{Andrews} and \cite{AndrewsLangfordMcCoy}) finally allow the construction of a complete Morse homology theory, at least in the case where the ambient manifold, $M$, is a flat, $(d+1)$-dimensional torus, and where certain further restrictions are also imposed on elements of $\Cal{E}(M)$.
\par
In \cite{SmiEMCFI} and \cite{SmiEMCFII}, we initiated a programme for the study of the Morse homology of the functional, $\Cal{F}$. The objective of the current paper is partly to improve on these results, but mainly to complete the final theoretical step required to complete the construction (c.f. Theorem \procref{ConcentrationParabolicCase}, below). The remaining work involved, which, though long and technical, is essentially formal, will be completed in our forthcoming paper, \cite{MacSmiMT}. There we will prove the following result, which we already state here in order to clearly illustrate our motivations.
\proclaim{Theorem \nextprocno, {\bf In preparation.}}
\noindent If $d\geq 2$ and if $T^{d+1}$ is a $(d+1)$-dimensional torus, then, for generic, smooth functions, $F:T^{d+1}\rightarrow]0,\infty[$, such that
$$
\msup_{x\in T^{d+1},\|\xi\|=1}\left|D^2F(x)(\xi,\xi)\right| < (3-2\sqrt{2})\minf_{x\in T^{d+1}}F(x)^3,
$$
there exist at least $2^{d+1}$ distinct Alexandrov-embeddings, $[e]\in\Cal{E}(T^{d+1})$, of mean curvature prescribed at every point by $F$.
\endproclaim
\proclabel{MainExistenceResult}
\remark Here a $(d+1)$-dimensional torus is taken to be {\sl any} quotient of $\Bbb{R}^{d+1}$ by a cocompact lattice.
\medskip
\remark Intriguingly, Theorem \procref{Separation}, below, would suggest that the Morse homology itself exhibits some sort of bifurcation behaviour as the forcing term moves beyond being subcritical. However, a deeper understanding of this lies far beyond the scope of the current paper.
\medskip
\remark An analogous result should also hold, with suitably modified conditions on $F$, when the ambient space is $2$-dimensional, and also when it is a compact, hyperbolic manifold.
\newsubhead{Morse homology}[MorseHomology]
In the infinite-dimensional setting, Morse homology algebraically encodes the relationship between the moduli spaces of solutions of certain non-linear para\-bolic operators, and the moduli spaces of solutions of those elliptic operators which correspond to their stationary states. However, in order to gain some intuition, it is worthwhile to first review the finite-dimensional theory (c.f. \cite{Schwarz} for a complete and thorough exposition). Thus, a smooth function, $f:M\rightarrow\Bbb{R}$, will be said to be of {\bf Morse} type whenever every one of its critical points is non-degenerate. It will then be said to be of Morse-Smale type whenever, in addition, every one of its complete gradient flows is non-degenerate in the following sense. Recall that a complete gradient flow, $\gamma:\Bbb{R}\rightarrow M$, of $f$ is, by definition, a zero of the non-linear differential operator
$$
\gamma\mapsto D_t\gamma + \nabla f(\gamma(t)).\eqnum{\nexteqnno[GradientFlowOperator]}
$$
The linearisation of this operator about $\gamma$, which maps sections of the pull-back bundle, $\gamma^*TM$, to other sections of the same bundle, is given by
$$
\xi\mapsto \nabla_{\partial_t}\xi + \opHess(f)(\gamma(t))\xi.\eqnum{\nexteqnno[LinearisedGradientFlowOperator]}
$$
The Morse property of $f$ ensures that this operator is always of Fredholm type, and the function, $f$, will then be of {\bf Morse-Smale} type whenever this operator is surjective for all $\gamma$. There is no shortage of functions with this property. Indeed, standard transversality results show that the set of all such functions is generic, that is, of the second category in the sense of Baire. This means that it contains the intersection of a countable family of open, dense subsets of $C^\infty(M)$, so that, in particular, by the Baire category theorem, it is dense.
\par
By the Morse property of $f$, every one of its critical points is isolated, so that the set, $Z$, of all its critical points is finite. In particular, for all $k$, the subset, $Z_k$, of $Z$, defined to be the set of critical points of Morse index $k$, is also finite, where the {\bf Morse index} of a critical point is here defined to be the sum of the multiplicities of all strictly negative eigenvalues of $\opHess(f)$ at that point. For all $k$, the $k$'th order chain group of the Morse homology of $(M,f)$ will then be defined by
$$
C_k(M,f) := \Bbb{Z}_2[Z_k].
$$
Equivalently, $C_k(M,f)$ will be the $\Bbb{Z}_2$-module of all $\Bbb{Z}_2$-valued functions over $Z_k$. In particular, the sum of the dimensions of the chain groups is equal to the cardinality of the solution set, $Z$.
\par
We will denote by $W$ the space of all complete gradient flows of $f$, furnished with the topology of $C^k_\oploc$-convergence for all $k$. By compactness and the Morse property of $f$, every complete gradient flow, $\gamma$, has well-defined end-points, $p$ and $q$, in the sense that
$$\eqalign{
\mlim_{t\rightarrow-\infty}\gamma(t) &= p,\ \text{and}\cr
\mlim_{t\rightarrow+\infty}\gamma(t) &= q.\cr}
$$
Furthermore, these end-points are always critical points of $f$, so that
$$
W = \munion_{(p,q)\in Z\times Z} W_{p,q},
$$
where $W_{p,q}$ here denotes the space of all complete gradient flows of $f$ starting at $p$ and ending at $q$. The Morse-Smale property of $f$ ensures that, for all $p,q$, $W_{p,q}$ in fact has the structure of a smooth manifold of dimension equal to the difference between the respective Morse indices of $p$ and $q$. Of particular interest is the case where this difference is equal to $1$, and where the space $W_{p,q}$ is therefore a one-dimensional manifold. Indeed, since the differential operator \eqnref{GradientFlowOperator} defining elements of $W_{p,q}$ is homogeneous in time, the additive group, $\Bbb{R}$, acts on this space by translation of the time variable, yielding a quotient space which happens to be a compact, zero-dimensional manifold: that is, a {\sl finite set of points}. The cardinality of this set will then be used to define the boundary operator of the chain complex of $(M,f)$. Indeed, for all $p\in Z_k$,
$$
\partial_k p := \sum_{q\in Z_{k-1}}\left[\# W_{p,q}/\Bbb{R}\right]q.
$$
\par
The first main theorem of Morse homology theory states that, for all $k$, the composition, $\partial_{k-1}\circ\partial_k$, vanishes. It follows that the chain complex, $(C_*(M),\partial_*)$, has a well-defined homology,
$$
H_*(M,f) := \opKer(\partial_*)/\opIm(\partial_{*+1}),
$$
and this will be called the {\bf Morse homology} of $(M,f)$. The second and third main theorems of Morse homology theory then state respectively that the Morse homology is, up to isomorphism, independent of the Morse-Smale function used, and, furthermore, that it is in fact isomorphic to the singular homology of the ambient space, $M$. Significantly, since the sum of the dimensions of the homology groups yields a lower bound for the sum of the dimensions of the chain groups, it also yields a lower bound for the cardinality of the solution set, $Z$. In this manner, it is shown that the number of critical points of a generic function, $f$, is bounded below by the sum of the Betti numbers of the ambient manifold, which is one of the main results of Morse theory.
\newsubhead{Eternal forced mean curvature flows}[EternalForcedMeanCurvatureFlows]
Given a riemannian manifold, $M$, and a smooth, positive function, $F:M\rightarrow]0,\infty[$, an {\bf eternal forced mean curvature flow} with forcing term $F$ will be a strongly smooth\footnote*{c.f. Appendix \headref{WeaklySmoothManifolds} for terminology.} curve $[e_t]:\Bbb{R}\rightarrow\Cal{E}(M)$ such that
$$
\langle\partial_t e_t, N_t\rangle + H_t - F\circ e_t = 0,\eqnum{\nexteqnno[ForcedMCFEquation]}
$$
where $N_t$ is the outward-pointing, unit, normal vector field of the embedding, $e_t$, and $H_t$ is its mean curvature. The eternal forced mean curvature flows with forcing term $F$ are precisely the $L^2$-gradient flows of the ``area-minus-volume'' functional, $\Cal{F}$, introduced in Section \subheadref{PrescribedCurvatureProblems}. Thus, bearing in mind the discussion of the preceding section, a Morse homology theory for the pair $(\Cal{E}(M),\Cal{F})$ should follow once the appropriate properties of eternal forced mean curvature flows have been established. The objective of the current paper will be to obtain precisely these properties in the case where $M$ is a $(d+1)$-dimensional torus, $T^{d+1}$.\footnote{${}^\dagger$}{Throughout this text, a $(d+1)$-dimensional torus will be taken to be {\sl any} quotient of $\Bbb{R}^{d+1}$ by a cocompact lattice.} We will first refine what has already been proven in \cite{SmiEMCFI} and \cite{SmiEMCFII}, before stating the new results which complete the theoretical side of the construction. Observe, in particular, that it will be sufficient to obtain results for solutions of the {\sl parabolic problem}, that is, for eternal forced mean curvature flows, since solutions of the {\sl elliptic problem}, that is, Alexandrov-embedded hyperspheres of prescribed curvature, arise as special cases of the former, being merely those forced mean curvature flows which are constant.
\par
It will first be necessary to introduce geometric restrictions. Thus, for all $\lambda\geq 1$, the subspace $\Cal{E}_{\lambda}(T^{d+1})$ of $\Cal{E}(T^{d+1})$ will be identified as follows. First, an embedding, $[e]$, will be said to be {\bf locally strictly convex} (or LSC) whenever every one of its principal curvatures is strictly positive. An LSC Alexandrov-embedding, $[e]$, will then be said to be {\bf pointwise $\lambda$-pinched} whenever it has the property that, for every point $x\in S^d$,
$$
\kappa_d(x) \leq \lambda\kappa_1(x),
$$
where $0<\kappa_1(x)\leq\kappa_d(x)$ are, respectively, the least and greatest principal curvatures of the Alexandrov-embedding, $e$, at the point $x$. Next recall that any LSC Alex\-endrov-embedding in $\Bbb{R}^{d+1}$ is actually the boundary of some open, convex set (c.f. \cite{Alexander}). A pointwise $\lambda$-pinched Alex\-androv-embedding, $[e]$, in $\Bbb{R}^{d+1}$ will then be said to be $\lambda$-non-collapsed whenever it has the property that for all $x\in S^d$, the Euclidean hypersphere of curvature $\lambda\kappa_1(x)$ which is an interior tangent to $[e]$ at the point $x$ is entirely contained within the closed set bounded by $[e]$.\footnote*{We alert the reader to the fact that our notion of non-collapsedness compares the curvature of the Euclidean hypersphere to the {\sl least} principal curvature of the Alexandrov-embedding. Although this may be at first surprising, it makes goods sense, since pointwise $\lambda$-pinching ensures that the curvature of the Euclidean hypersphere is actually no less than the greatest principal curvature of the Alexandrov-embedding.} A pointwise $\lambda$-pinched Alexandrov-embedding, $[e]$, in $T^{d+1}$ is then said to be {\bf $\lambda$-non-collapsed} whenever its lift to $\Bbb{R}^{d+1}$ has this property. Finally, for all $\lambda\geq 1$, the subspace $\Cal{E}_\lambda(T^{d+1})$ of $\Cal{E}(T^{d+1})$ will be defined to be the set of all Alexandrov-embeddings, $[e]$, which are LSC, pointwise $\lambda$-pinched, and $\lambda$-non-collapsed.
\par
Consider now the rational function $\phi:[1,\infty[\rightarrow[0,\infty[$ given by
$$
\phi(t) := \frac{(t-1)}{t(t+1)}.\eqnum{\nexteqnno[ControllingFunction]}
$$
Observe that this function tends to $0$ at $1$ and $+\infty$ and has a unique maximum, equal to $(3-2\sqrt{2})$, at the point $(1+\sqrt{2})$. There therefore exist two smooth inverses,
$$\eqalign{
\lambda:[0,(3-2\sqrt{2})[&\rightarrow[1,1+\sqrt{2}[,\ \text{and}\cr
\Lambda:[0,(3-2\sqrt{2})[&\rightarrow]1+\sqrt{2},+\infty].\cr}
$$
The forcing term, $F$, will be said to be {\bf sub-critical} whenever it is strictly positive, and
$$
\msup_{x\in T^{d+1},\|\xi\|=1}\left|D^2(F)(x)(\xi,\xi)\right| \leq (3-2\sqrt{2})\minf_{x\in T^{d+1}}F(x)^3.\eqnum{\nexteqnno[Subcritical]}
$$
In this case, the constants $\lambda_F<\Lambda_F$ will be defined by
$$\eqalign{
\lambda_F &:= \lambda\left(\msup_{x\in T^{d+1},\|\xi\|=1}\left|D^2(F)(x)(\xi,\xi)\right|/\minf_{x\in T^{d+1}}F(x)^3\right),\ \text{and}\cr
\Lambda_F &:= \Lambda\left(\msup_{x\in T^{d+1},\|\xi\|=1}\left|D^2(F)(x)(\xi,\xi)\right|/\minf_{x\in T^{d+1}}F(x)^3\right),\cr}\eqnum{\nexteqnno[LittleLambdaAndLambda]}
$$
and for subcritical $F$, the subspace $\Cal{E}_F(T^{d+1})$ will be defined by
$$
\Cal{E}_F(T^{d+1}) := \Cal{E}_{\Lambda_F}(T^{d+1}).\eqnum{\nexteqnno[AdmissableElements]}
$$
Embeddings in $\Cal{E}_F(T^{d+1})$ will be called {\bf admissable}. Likewise, a forced mean curvature flow, $[e_t]$, with forcing term, $F$, will be said to be {\bf admissable} whenever $[e_t]$ is admissable for all $t$. Morse homology will be constructed for the pair $(\Cal{E}_F(T^{d+1}),\Cal{F})$.
\par
Significantly, the space, $\Cal{E}_F(T^{d+1})$, of admissable Alexandrov-embeddings is actually a {\sl closed} subset of $\Cal{E}(T^{d+1})$. However, it is of fundamental importance in any differential topology theory that the space of admissable elements be {\sl open}. Indeed, otherwise, the perturbative stages of the construction would cease to function. This problem of openness is typically addressed indirectly in the statement of the compactness result, which is the case in our earlier work, \cite{SmiEMCFI}, where we prove compactness for families of eternal forced mean curvature flows, and where the rather restrictive conditions imposed actually serve to ensure that all limits remain within a given open set. The greater generality of the present setting is obtained via the following result, where the problem of openness is addressed via an adaptation of Ben Andrews' non-collapsing theorem. Indeed, even though the subset, $\Cal{E}_F(T^{d+1})$, is not open, for subcritical $F$, the Morse homology of $(\Cal{E}_F(T^{d+1}),\Cal{F})$ lies strictly in its interior, and is, in particular, separated from the rest of the Morse homology of $(\Cal{E}(T^{d+1}),\Cal{F})$. An eternal forced mean curvature flow $[e_t]:\Bbb{R}\rightarrow\Cal{E}(T^{d+1})$ will be said to be of {\bf bounded type} whenever
$$
\msup_{t}\opDiam([e_t]) < \infty.
$$
In Section \headref{Admissability}, we prove
\proclaim{Theorem \nextprocno, {\bf Separation.}}
\noindent If $[e_t]:\Bbb{R}\rightarrow\Cal{E}_F(T^{d+1})$ is an eternal forced mean curvature flow of bounded type with sub-critical forcing term $F$, and if $[e_t]$ is pointwise $\Lambda_F$-pinched and $\Lambda_F$-non-collapsed for all $t$, then $[e_t]$ is pointwise $\lambda_F$-pinched and $\lambda_F$-non-collapsed for all $t$.
\endproclaim
\proclabel{Separation}
Theorem \procref{Separation}, justifies the context of all that follows. First, the following compactness result is obtained via a straightforward blow-up argument (c.f. Section \subheadref{Compactness}).
\proclaim{Theorem \nextprocno, {\bf Compactness.}}
\noindent Fix $d\geq 2$, and let $(F_m)$ be a sequence of smooth, positive functions over $\Bbb{R}^{d+1}$ converging in the $C^k_\oploc$-sense for all $k$ to the smooth, positive function $F_\infty$. Suppose, furthermore, that
$$
0 < \minf_m F_{m,-} \leq \msup_m F_{m,+} < \infty.
$$
For all $m$, let $[e_{m,t}]:\Bbb{R}\rightarrow\Cal{E}(\Bbb{R}^{d+1})$ be an eternal forced mean curvature flow of bounded type with forcing term $F_m$. Suppose, furthermore, that there exists $\lambda\geq 1$ such that, for all $m$, and for all $t$, $[e_{m,t}]$ is pointwise $\lambda$-pinched and $\lambda$-non-collapsed. If there exists a compact subset $K\subseteq\Bbb{R}^{d+1}$ such that $[e_{m,0}]\minter K\neq\emptyset$ for all $m$, then there exists an eternal forced mean curvature flow, $[e_{\infty,t}]$, towards which the sequence $([e_{m,t}])$ subconverges in the $C^k_\oploc$-sense for all $k$. In particular, $[e_{\infty,t}]$ is also of bounded type and, for all $t$, $[e_{\infty,t}]$ is pointwise $\lambda$-pinched and $\lambda$-non-collapsed.
\endproclaim
\proclabel{Compactness}
\noindent As in the finite-dimensional case, the functional, $\Cal{F}$, will be said to be of {\bf Morse type} whenever every one of its admissable critical points is non-degenerate. Recalling that Theorem \procref{Compactness} also applies to sequences of critical points of $\Cal{F}$, straightforward differential-topological arguments now show (c.f., for example, \cite{RosSmiDT} and \cite{WhiteI}) that the functional, $\Cal{F}$, has this property for generic $F$. In this case, it readily follows from Theorem \procref{Compactness}, again, that eternal forced mean curvature flows always have well-defined end-points.
\proclaim{Theorem \nextprocno, {\bf End-points.}}
\noindent If $\Cal{F}$ is of Morse type, and if $[e_t]:\Bbb{R}\rightarrow\Cal{E}_F(T^{d+1})$ is an admissable eternal forced mean curvature flow of bounded type, then there exist admissable embeddings, $[e_\pm]\in\Cal{E}_F(T^{d+1})$, of mean curvature prescribed by $F$ such that
$$\eqalign{
\mlim_{t\rightarrow -\infty}[e_t] &= [e_-],\ \text{and}\cr
\mlim_{t\rightarrow +\infty}[e_t] &= [e_+].\cr}
$$
\endproclaim
\proclabel{EndPoints}
\par
However, the main result of this paper actually concerns the asymptotic behaviour of the Morse homology of $(\Cal{E}_F(T^{d+1}),\Cal{F})$ as $F$ tends to infinity. Indeed, it is this that will yield an explicit formula for the Morse homology. Consider therefore a smooth function, $f:T^{d+1}\rightarrow\Bbb{R}$, of Morse-Smale type. For all sufficiently small $\kappa>0$, consider the forcing term
$$
F_\kappa := \frac{1}{\kappa}(1 - \kappa^2f),
$$
let $\Cal{F}_\kappa$ denote its ``area-minus-volume'' functional, and let $\Cal{E}_\kappa(T^{d+1})$ denote the space of all Alexandrov-embedded hyperspheres which are admissable for $F_\kappa$, that is
$$
\Cal{E}_\kappa(T^{d+1}) := \Cal{E}_{F_\kappa}(T^{d+1}).
$$
In \cite{SmiEC}, building on the work, \cite{Ye}, of Ye, we show (c.f. also \cite{PacardXu}),
\proclaim{Theorem \nextprocno, {\bf Concentration: elliptic case.}}
\noindent For sufficiently small $\kappa>0$, there exists a strongly smooth map $\Phi:T^{d+1}\rightarrow\Cal{E}_\kappa(T^{d+1})$ such that for any critical point, $p$, of $f$ of Morse index $k$, the point $\Phi(p)$ is a non-degenerate critical point of $\Cal{F}_\kappa$ of Morse index $(k+1)$. Furthermore, $\Cal{F}_\kappa$ has no other critical points in $\Cal{E}_\kappa(T^{d+1})$.
\endproclaim
\proclabel{ConcentrationEllipticCase}
The Morse homology of $(\Cal{E}_\kappa(T^{d+1}),\Cal{F}_\kappa)$ is explicitly determined once the corresponding asymptotic behaviour for eternal forced mean curvature flows has been proven. Thus, given any two critical points, $[e]$ and $[f]$, of $\Cal{F}_\kappa$, $\Cal{W}_{[e],[f]}$ will denote the space of admissable forced mean curvature flows of bounded type with forcing term $F_\kappa$, starting at $[e]$ and ending at $[f]$. Next, an eternal forced mean curvature flow, $[e_t]:\Bbb{R}\rightarrow\Cal{E}(T^{d+1})$, will be said to be {\bf non-degenerate} whenever its linearised mean curvature flow operator is surjective. If every element of $\Cal{W}_{[e],[f]}$ is non-degenerate, and if the functional, $\Cal{F}_\kappa$, is of Morse type, then it follows by standard techniques of Fredholm theory that $\Cal{W}_{[e],[f]}$ has the structure of a smooth, finite-dimensional manifold. In Sections \headref{SingularPerturbation} and \headref{Concentration}, making use of the compactness result of Theorem \procref{Compactness}, together with the existence of well-defined end-points established in Theorem \procref{EndPoints}, we complement Theorem \procref{ConcentrationEllipticCase} by showing
\proclaim{Theorem \nextprocno, {\bf Concentration: parabolic case.}}
\noindent For sufficiently small $\kappa>0$, and for every pair of critical points, $p$ and $q$, of $f$ such that $\opIndex(p)-\opIndex(q)=1$, every element of $\Cal{W}_{\Phi(p),\Phi(q)}$ is non-degenerate, and there exists a canonical diffeomorphism
$$
\hat{\Phi}:W_{p,q}\rightarrow\Cal{W}_{\Phi(p),\Phi(q)}.
$$
In particular, there are no other admissable eternal forced mean curvature flows of bounded type with forcing term $F_\kappa$, starting at $\Phi(p)$ and ending at $\Phi(q)$.
\endproclaim
\proclabel{ConcentrationParabolicCase}
\remark This result follows directly from Theorem \procref{SingularPerturbations} and Lemma \procref{FlowsAreAsymptoticToFormalSolution}, below. In actual fact, in \cite{SmiEMCFII}, the perturbative part of this result has been shown in the slightly different context of forced mean curvature flows with constant forcing term inside general riemannian manifolds.
\newsubhead{Constructing the Morse homology}[ConstructingTheMorseHomology]
It remains to sketch how the results of Section \subheadref{EternalForcedMeanCurvatureFlows} serve to construct the Morse homology of $(\Cal{E}_F(T^{d+1}),\Cal{F})$. First, for a subcritical forcing term, $F$, the {\bf elliptic solution space}, $\Cal{Z}$, is defined to be the set of all critical points of $\Cal{F}$; that is, the set of all admissable, Alexandrov-embedded hyperspheres of mean curvature prescribed by $F$. As indicated above, it follows from the compactness result of Theorem \procref{Compactness} that, for generic $F$, the ``area-minus-volume'' functional, $\Cal{F}$, is of Morse type in the sense that every one of its critical points is non-degenerate. In particular, in this case, $\Cal{Z}$ consists only of isolated elements, and, by the compactness result of Theorem \procref{Compactness}, again, is finite.
\par
It follows from the spectral theory of second-order, elliptic operators that every element, $[e]$, of $\Cal{Z}$ has finite Morse index, defined here to be the sum of the multiplicities of all the strictly negative eigenvalues of the Jacobi operator of $\Cal{F}$ at $e$, (c.f. \cite{RosSmiDT} and \cite{WhiteI}). For all integer $k$, the finite set, $\Cal{Z}_k$, is then defined to be the set of all critical points of $\Cal{F}$ of Morse index $k$, and, the $k$'th order chain group of $(\Cal{E}_F(T^{d+1}),\Cal{F})$ is defined by
$$
\Cal{C}_k(\Cal{E}_F(T^{d+1}),\Cal{F}) := \Bbb{Z}_2[\Cal{Z}_k],
$$
so that, in particular, the sum of the dimensions of the chain groups is equal to the cardinality of the elliptic solution space, $\Cal{Z}$.
\par
The {\bf parabolic solution space}, $\Cal{W}$, is defined to be the space of all admissable, eternal forced mean curvature flows of bounded type and with forcing term $F$ furnished with the topology of $C^k_\oploc$-convergence for all $k$. By the existence of well-defined end-points established in Theorem \procref{EndPoints}, this space decomposes as
$$
\Cal{W} = \munion_{[e],[f]\in\Cal{Z}\times\Cal{Z}}\Cal{W}_{[e],[f]},
$$
where $\Cal{W}_{[e],[f]}$ here denotes the subspace of $\Cal{W}$ consisting of those flows which start at $[e]$ and which end at $[f]$. Upon perturbing $\Cal{F}$ in a straightforward manner, we may suppose that every element of $\Cal{W}_{[e],[f]}$ is non-degenerate, so that, by the inverse function theorem, the space, $\Cal{W}_{[e],[f]}$, carries the structure of a smooth manifold of finite-dimension, which, by the Atiyah-Patodi-Singer Index Theorem (c.f. \cite{RobbinSalamon}), is equal to the difference between the respective Morse indices of $[e]$ and $[f]$. In particular, when this difference is equal to $1$, $\Cal{W}_{[e],[f]}$ is $1$-dimensional, and since the equation, \eqnref{ForcedMCFEquation}, defining elements of $\Cal{W}$ is homogeneous in time, the additive group, $\Bbb{R}$, acts on this space by translation of the time variable, yielding a quotient space which, as before, turns out to be a compact, $0$-dimensional manifold: that is, a finite set of points. The cardinality of this set is used to define the boundary operator of the chain complex of $(\Cal{E}_F(T^{d+1}),\Cal{F})$. Indeed, for all $[e]\in\Cal{Z}_k$,
$$
\partial_l [e] := \sum_{[f]\in\Cal{C}_{l-1}}\left[\#\Cal{W}_{[e],[f]}/\Bbb{R}\right][f].
$$
\par
The key result of Morse homology theory is that the square of the boundary operator vanishes. The two main components of this result are a compactness theorem modulo broken trajectories and, conversely, a glueing theorem for broken trajectories. However, a straightforward combinatorial argument, outlined in \cite{SmiEMCFI}, shows how compactness modulo broken trajectories readily follows from the compactness result already obtained in Theorem \procref{Compactness}. The remaining steps, including the glueing theorem, then follow by general arguments valid for large families of suitably regular parabolic operators. It is therefore a straightforward, though highly technical, matter to prove the fundamental relation
$$
\partial_{l-1}\circ\partial_l = 0,
$$
for all $l$. It follows that the chain complex, $\Cal{C}_*(\Cal{E}_F(T^{d+1}),\Cal{F})$, has a well defined homology,
$$
\Cal{H}_*(\Cal{E}_F(T^{d+1}),\Cal{F}) := \opKer(\partial_*)/\opIm(\partial_{*+1}),
$$
which we call the Morse homology of $(\Cal{E}_F(T^{d+1}),\Cal{F})$.
\par
Similar arguments then show that, up to isomorphism, the Morse homology is independent of the generic, sub-critical function, $F$, used, and it is then explicitely calculated by considering the case where $F=F_\kappa$, for some Morse-Smale function, $f$, and for suitably small $\kappa$. Indeed, by Theorem \procref{ConcentrationEllipticCase}, for all $l$, and for all sufficiently small $\kappa$, there exists a canonical isomorphism
$$
\Phi_l:C_l(X,f)\rightarrow\Cal{C}_{l+1}(\Cal{E}_\kappa(T^{d+1}),\Cal{F}_\kappa).
$$
Furthermore, by Theorem \procref{ConcentrationParabolicCase}, upon reducing $\kappa$ further if necessary,
$$
\partial_{l+1}\circ\Phi_l = \Phi_{l-1}\circ\partial_l,
$$
so that $\Phi_l$ quotients down to another isomorphism
$$
\Phi_l:H_l(X,f)\rightarrow\Cal{H}_{l+1}(\Cal{E}_\kappa(T^{d+1}),\Cal{F}_\kappa).
$$
Since $H_*(X,f)$ is known to be isomorphic to the singular homology of the torus, it then follows that, for generic, subcritical $F$,
$$
H_l(\Cal{E}_F(T^{d+1}),\Cal{F})
=\left\{\matrix\Bbb{Z}_2^{\left({}^{d+1}_{l-1}\right)}\hfill&\ \text{if}\ 1\leq l\leq d+2\hfill\cr
0\hfill&\ \text{otherwise.}\hfill\cr\endmatrix\right.
$$
Finally, since the sum of the dimensions of the homology groups yields a lower bound for the sum of the dimensions of the chain groups, it also yields a lower bound for the cardinality of the elliptic solution space, $\Cal{Z}$, so that, for a generic, subcritical forcing term, $F$,
$$
\#\Cal{Z} \geq \sum_{l=1}^{d+2}\left({}^{d+1}_{l-1}\right) = 2^{d+1},
$$
thus proving Theorem \procref{MainExistenceResult}.
\newsubhead{Acknowledgements}[Acknowledgements]
The author is grateful to Matthias Schwarz for having drawn his attention to the possibility of applying Morse homology theory to curvature flows. He is also grateful to Harold Rosenberg for having suggested the application of differential-topological techniques to the study of hypersurface problems. Finally, this project has been rather long in gestation, and many of the ideas were conceived whilst the author was benefitting from  a
Marie Curie postdoctoral fellowship in the Centre de Recerca Matem\`atica, Barcelona.
\newhead{Compactness and end-points}[CompactnessAndEndPoints]
\newsubhead{Geometric bounds}[GeometricBounds]
Given a riemannian manifold, $M$, and a smooth, positive function $F:M\rightarrow]0,\infty[$, denote
$$\eqalign{
F_- &:= \minf_{x\in M}F(x),\ \text{and}\cr
F_+ &:= \msup_{x\in M}F(x).\cr}
$$
Consider now the case where the ambient space is $(d+1)$-dimensional Euclidean space.
\proclaim{Lemma \nextprocno}
\noindent If $[e_t]:\Bbb{R}\rightarrow\Cal{E}(\Bbb{R}^{d+1})$ is an eternal forced mean curvature flow of bounded type with forcing term $F$, then, for all $t$, $[e_t]$ contains no Euclidean hypersphere of mean curvature less than $F_-$.
\endproclaim
\proclabel{DoesNotContainBigSpheres}
\proof Suppose the contrary. Without loss of generality, $[e_0]$ contains the Euclidean hypersphere of mean curvature $a<F_-$ about $0$. Let $[f_t]:[0,\infty[\rightarrow\Cal{E}(\Bbb{R}^{d+1})$ be the forced mean curvature flow with constant forcing term, $F_-$, starting at this hypersphere. Explicitly,
$$
f(t,x) =  \psi^{-1}\left(F_-^2t + \psi\left(\frac{1}{a}\right)\right)x,
$$
where $\psi:]1/F_-,\infty[\rightarrow\Bbb{R}$ is given by
$$
\psi(r) := (rF_- - 1) + \opLog(rF_- - 1).
$$
Since $[f_0]$ is contained in $[e_0]$, and since $F_-\leq F$, it follows by the geometric maximum principle that $[f_t]$ is contained in $[e_t]$ for all $t$. However, this is not possible, since $[e_t]$ is of bounded type, but $[f_t]$ is unbounded, and the result follows.\qed
\proclaim{Lemma \nextprocno}
\noindent If $[e_t]:\Bbb{R}\rightarrow\Cal{E}(\Bbb{R}^{d+1})$ is an eternal forced mean curvature flow with forcing term, $F$, then, for all $t$, $[e_t]$ is contained within no Euclidean hypersphere of mean curvature greater than $F_+$.
\endproclaim
\proclabel{IsNotContainedInSmallSpheres}
\proof Suppose the contrary. Without loss of generality, $[e_0]$ is contained within the Euclidean hypersphere of curvature $a>F_+$ about $0$. Let $[f_t]:[0,R[\rightarrow\Cal{E}(\Bbb{R}^{d+1})$ be the forced mean curvature flow with constant forcing term, $F_+$, starting at this hypersphere. Explicitely,
$$
f(t,x) = \psi^{-1}\left(F_+^2t + \psi\left(\frac{1}{a}\right)\right)x,
$$
where $\psi:]0,1/F_+[\rightarrow\Bbb{R}$ is given by
$$
\psi(r) := (rF_+ - 1) + \opLog(1 - rF_+).
$$
Since $[f_0]$ contains $[e_0]$ and since $F_+\geq F$, it follows by the geometric maximum principle that $[f_t]$ contains $[e_t]$ for all $t$. However, this is not possible, since $[e_t]$ exists for all time, but $[f_t]$ extinguishes in finite time, and the result follows.\qed
\medskip
Given a constant, $c>0$, an LSC embedding, $[e]$ is said to be {\bf $c$-convex} whenever its least principal curvature is at every point bounded below by $c$.
\proclaim{Lemma \nextprocno}
\noindent Let $[e_t]:\Bbb{R}\rightarrow\Cal{E}(\Bbb{R}^{d+1})$ be an eternal forced mean curvature flow of bounded type with forcing term $F$. If $[e_t]$ is $\lambda$-non-collapsed, then it is also $c$-convex, where
$$
c = \frac{F_-}{\lambda}.
$$
\endproclaim
\proclabel{CConvexity}
\proof Indeed, fix $t\in\Bbb{R}$ and let $x\in S^d$ minimise $\kappa_1$. Since $[e_t]$ is $\lambda$-non-collapsed, the hypersphere of curvature $\lambda\kappa_1(x)$ which is tangent to $[e_t]$ at $x$ is contained within $[e_t]$. Thus, by Lemma \procref{DoesNotContainBigSpheres}, $\lambda\kappa_1(x)>F_-$, and the result follows.\qed
\medskip
\noindent In particular, the Bonnet-Meyer Theorem yields explicit bounds on the diameter of each $[e_t]$.
\proclaim{Corollary \nextprocno}
\noindent With the same hypotheses as in Lemma \procref{CConvexity}, for all $t$,
$$
\opDiam([e_t]) \leq \sqrt{\frac{\lambda}{F_-}}\pi.
$$
\endproclaim
\proclabel{BoundedDiameter}
\newsubhead{Compactness}[Compactness]
Consider a sequence, $(F_m)$, of smooth, positive functions over $\Bbb{R}^{d+1}$ converging in the $C^k_\oploc$-sense for all $k$ to the smooth, positive function, $F_\infty$. Suppose, furthermore, that
$$
0 < \minf_m F_{m,-} \leq \msup_m F_{m,+} < \infty,
$$
where $F_{m,-}$ and $F_{m,+}$ are defined as in the preceding section. For all $m$, let $[e_{m,t}]:\Bbb{R}\rightarrow\Cal{E}(\Bbb{R}^{d+1})$ be an eternal forced mean curvature flow of bounded type with forcing term $F_m$, and suppose that there exists a compact subset, $K$, of $\Bbb{R}^{d+1}$ such that
$$
[e_{m,0}] \minter K \neq \emptyset
$$
for all $m$. Recall first the classical
\proclaim{Lemma \nextprocno, {\bf Quasi-maximum lemma.}}
\noindent Let $X$ be a complete metric space. Let $\phi:X\rightarrow]0,\infty[$ be a positive function. For all $x_0\in X$, for all $C>1$, and for all $A,\alpha>0$, there exists $x$ in $X$ such that $\phi(x)\geq \phi(x_0)$, and
$$
\phi(y) \leq C\phi(x)
$$
for all $y\in B_{A\phi(x)^{-\alpha}}(x)$, where $B_{A\phi(x)^{-\alpha}}(x)$ here denotes the ball of radius $A\phi(x)^{-\alpha}$ about the point, $x$.
\endproclaim
\proclabel{QuasiMaximumLemma}
\proof Indeed, otherwise, by induction, there exists a sequence, $(x_m)$, of points in $X$ such that, for all $m$, $x_{m+1}\in B_{A\phi(x)^{-\alpha}}(x_m)$, and $\phi(x_{m+1})>C\phi(x_m)$. By induction again, for all $m$, $\phi(x_m)>C^m\phi(x)$, and so $d(x_{m+1},x_m)<A\phi(x_0)^{-\alpha}C^{-m\alpha}$. Since the series $\sum_{m=0}^\infty C^{-m\alpha}$ converges, it follows that $(x_m)$ is a Cauchy sequence in $X$. By completeness, there therefore exists $x_\infty$ in $X$ towards which $(x_m)$ converges. However, this is absurd, since the sequence $(\phi(x_m))$ diverges, but $\phi$ is continuous. The result follows.\qed
\proclaim{Lemma \nextprocno}
\noindent If $d\geq 2$, and if, for some $\lambda\geq 1$, $[e_{m,t}]$ is pointwise $\lambda$-pinched and $\lambda$-non-collapsed for all $m$, then there exists $B>0$ such that, for all $m$, and for all $(t,x)$,
$$
\kappa_{m,d,t}(x)\leq B,
$$
where $\kappa_{m,d,t}$ here denotes the greatest principal curvature of the embedding $e_{m,t}$.
\endproclaim
\proclabel{BoundedCurvature}
\proof Suppose the contrary. For all $m$, let $K_m:\Bbb{R}\rightarrow\Bbb{R}$ be the function given by
$$
K_m(t) = \msup_{x\in S^d}\kappa_{m,t,d}(x),
$$
and consider a sequence $(t_m)$, chosen such that $(K_m(t_m))$ tends to $+\infty$. By the quasi-maximum lemma (Lemma \procref{QuasiMaximumLemma}), without loss of generality, for all $m$, and for all $\left|s-t_m\right|<K_m(t_m)^{-1/2}$,
$$
K_m(s) < 2K_m(t_m).
$$
For all $m$, let $x_m$ be a point of $S^d$ maximising $K_d$, and define $[\tilde{e}_{m,t}]$ by
$$
\tilde{e}_{m,t}(x) := K_m(t_m)\left(e_{m,(t-t_m)/K_m(t_m)}(x) - e_{m,t_m}(x_m)\right),
$$
so that, for all $m$, $[\tilde{e}_{m,t}]$ is an eternal forced mean curvature flow of bounded type with forcing term
$$
\tilde{F}_{m,t} := \frac{1}{K_m(t_m)}\left(F_m(x/K_m(t_m) + x_m)\right).
$$
Furthermore, for all $m$, $\tilde{e}_m(0,x_m)=0$, $\tilde{\kappa}_{m,d}(0,x_m)=1$, and
$$
\tilde{\kappa}_{m,d}(t,x) \leq 2
$$
for all $(t,x)\in[-K_m(t_m)^{1/2},K_m(t_m)^{1/2}]\times S^d$.
\par
For all $m$, let $X_m$ be the convex subset bounded by the embedded hypersphere $[\tilde{e}_{m,0}]$. Let $X_\infty$ be a convex subset towards which $(X_m)$ subconverges in the local Hausdorff sense. It follows by the above estimates and hypoelliptic regularity that the sequence, $(X_m)$, actually subconverges towards $X_\infty$ in the $C^k_\oploc$-sense for all $k$. In particular, $\partial X_\infty$ is convex, pointwise $\lambda$-pinched and $\lambda$-non-collapsed, and its greatest principal curvature at the point $x_\infty$ is equal to $1$. In particular, $X_\infty$, is strictly convex, since, otherwise, being $\lambda$-non-collapsed, it would coincide with a half-space and therefore have zero curvature at every boundary point, which is absurd. Since $d\geq 2$, it follows by the result, \cite{Hamilton}, of Hamilton that $X_\infty$ is compact, and there therefore exists $C>0$ such that, after extracting a subsequence,
$$
\opDiam([e_{m,t_m}]) \leq C/K_m(t_m)
$$
for all $m$. In particular, for sufficiently large $m$, $[e_{m,t_m}]$ is contained within a Euclidean hypersphere of arbitrarily large curvature. This is absurd by Lemma \procref{IsNotContainedInSmallSpheres}, and the result follows.\qed
\medskip
\noindent This proves the compactness result for eternal forced mean curvature flows.
\proclaim{Theorem \procref{Compactness}, {\bf Compactness}}
\noindent If $d\geq 2$, and if, for some $\lambda\geq 1$, $[e_{m,t}]$ is pointwise $\lambda$-pinched and $\lambda$-non-collapsed for all $m$, then there exists an eternal forced mean curvature flow, $[e_{\infty,t}]$, of bounded type with forcing term $F_\infty$ towards which $[e_{m,t}]$ subconverges in the $C^k_\oploc$-sense for all $k$. In particular, $[e_{\infty,t}]$ is also pointwise $\lambda$-pinched and $\lambda$-non-collapsed.
\endproclaim
\proof By standard arguments from the theory of quasi-linear parabolic PDEs, it suffices to show that the diameters and principal curvatures of all the spheres $([e_{m,t}])$ are uniformly bounded above. The result now follows by Corollary \procref{BoundedDiameter} and Lemma \procref{BoundedCurvature}.\qed
\newsubhead{End points}[EndPoints]
Consider now the case where the ambient manifold is a $(d+1)$-dimensional torus. Consider a smooth, positive function, $F:T^{d+1}\rightarrow]0,\infty[$ and admissible, eternal forced mean curvature flows in $T^{d+1}$ with forcing term $F$. First recall
\proclaim{Lemma \nextprocno, {\bf Unique continuation.}}
\noindent If $[e_t],[f_t]:\Bbb{R}\rightarrow\Cal{E}(\Bbb{T}^{d+1})$ are eternal forced mean curvature flows with forcing term $F$, and if $[f_t]=[e_{t+\delta}]$ for some $t$ and for some $\delta$, then $[f_t]=[e_{t+\delta}]$ for all $t$.
\endproclaim
\proclabel{UniqueContinuation}
\proof This follows by a standard argument using the result, \cite{AgmonNirenberg}, of Agmon \& Nirenberg (c.f. \cite{Weber} for details). Observe, in particular, that energy-boundedness is not a concern for us, since the flows considered here are all smooth.\qed
\medskip
\noindent Denote by $\Cal{F}([e_t])$ the set $\left\{\Cal{F}([e_t])\ |\ t\in\Bbb{R}\right\}$.
\proclaim{Lemma \nextprocno}
\noindent Let $[e_t]:\Bbb{R}\rightarrow\Cal{E}(T^{d+1})$ be an admissible, eternal forced mean curvature flow of bounded type with forcing term $F$. If $[e_t]$ is non-constant, then the set $\Cal{F}([e_t])$ is an open, bounded interval.
\endproclaim
\proclabel{EitherConstantOrOpen}
\proof Consider the smooth function, $\phi:\Bbb{R}\rightarrow\Bbb{R}$, given by $\phi(t)=\Cal{F}([e_t])$. In particular, $\Cal{F}([e_t])=\opIm(\phi)$. Since $\phi$ is continuous, $\opIm(\phi)$ is connected, and is therefore an interval. Since $[e_t]$ is an $L^2$-gradient flow of $\Cal{F}$, $\phi'(t)=0$ if and only if $[e_t]$ is a critical point of $\Cal{F}$. It follows by the unique continuation result of Lemma \procref{UniqueContinuation}, that if $\phi'$ vanishes at any point, then $[e_t]$ is the constant flow sending $\Bbb{R}$ to a critical point of $\Cal{F}$. Since this case is excluded by hypothesis, $\phi'$ never vanishes, and so $\opIm(\phi)$ is open.
\par
It remains to show that $\Cal{F}([e_t])$ is bounded. Suppose, first, that it is not bounded above. Then there exists a sequence $(t_m)$ (converging to $-\infty$) such that $(\Cal{F}([e_{t_m}]))$ tends to infinity. For all $m$, define
$$
\tilde{e}_m(t,x) := e_m(t-t_m,x),
$$
so that $[\tilde{e}_{m,t}]$ is also an admissable, eternal forced mean curvature flow of bounded type with forcing term $F$. By the compactness result of Theorem \procref{Compactness}, there exists an eternal forced mean curvature flow, $[\tilde{e}_{\infty,t}]$, also with forcing term $F$, towards which $[\tilde{e}_{m,t}]$ subconverges in the $C^k_\oploc$-sense for all $k$. However,
$$
\Cal{F}_\infty([\tilde{e}_{\infty,0}]) = \mlim_{m\rightarrow\infty}\Cal{F}_m([\tilde{e}_{m,0}])
=\mlim_{m\rightarrow\infty}\Cal{F}([e_{m,t_m}]) = +\infty,
$$
which is absurd, and it follows that $\Cal{F}([e_t])$ is bounded above. In the same manner, $\Cal{F}([e_t])$ is shown to be bounded below, and this completes the proof.\qed
\proclaim{Lemma \nextprocno}
\noindent Let $[e_t]:\Bbb{R}\rightarrow\Cal{E}(T^{d+1})$ be an admissible, eternal forced mean curvature flow of bounded type with forcing term $F$. For any sequence $(t_m)$ converging to $\pm\infty$, there exists an admissible critical point, $[e_\infty]$, of $\Cal{F}$ towards which the sequence $([e_{t_m}])$ subconverges.
\endproclaim
\proclabel{EndPointsExist}
\proof Observe that the sequence $(\Cal{F}([e_{t_m}]))$ converges to one of the two boundary points of $I:=\Cal{F}([e_t])$. For all $m$, define
$$
\tilde{e}_m(t,x) := e_m(t-t_m,x),
$$
so that $[\tilde{e}_m]$ is also an admissible, eternal forced mean curvature flow of bounded type with forcing term $F$. By the compactness result of Theorem \procref{Compactness}, there exists an admissible, eternal forced mean curvature flow, $[\tilde{e}_{\infty,t}]$, also with forcing term $F$, towards which $[\tilde{e}_m]$ subconverges in the $C^k_\oploc$-sense for all $k$. However, $\Cal{F}([\tilde{e}_{\infty,t}])\subseteq I$, but
$$
\Cal{F}([\tilde{e}_{\infty,0}]) = \mlim_{m\rightarrow\infty}\Cal{F}([\tilde{e}_{m,0}])
=\mlim_{m\rightarrow\infty}\Cal{F}([\tilde{e}_{m,t_m}])\in\partial I,
$$
so that $\Cal{F}([\tilde{e}_\infty])$ cannot be open. It follows by Lemma \procref{EitherConstantOrOpen}, that $[\tilde{e}_{\infty,t}]$ is a constant flow mapping $\Bbb{R}$ to a critical point of $\Cal{F}$. In particular, the sequence, $([e_{t_m}])$, subconverges to $[\tilde{e}_{\infty,0}]$, and the result follows.\qed
\medskip
The $A$- and $\Omega$- limit sets of the flow $[e_t]$ are defined by
$$\eqalign{
A([e_t]) &:= \minter_{t\in\Bbb{R}}\overline{\left\{[e_s]\ |\ s\leq t\right\}},\ \text{and}\cr
\Omega([e_t]) &:= \minter_{t\in\Bbb{R}}\overline{\left\{[e_s]\ |\ s\geq t\right\}}.\cr}
$$
By Lemma \procref{EndPointsExist}, these sets are both non-empty and consist only of admissible critical points of $\Cal{F}$. Furthermore, being intersections of nested families of connected sets, they are also both connected.
\par
As indicated in the introduction, since admissable critical points of $\Cal{F}$ are, in particular, constant eternal forced mean curvature flows, the compactness result of Theorem \procref{Compactness} also applies to these objects. Standard differential-topological arguments then show (c.f., for example, \cite{RosSmiDT}) that, for a generic choice of the density, $F$, the functional, $\Cal{F}$, is of {\bf Morse type} in the sense that all of its admissible critical points are non-degenerate. In particular, in this case, by connectedness, the $A$- and $\Omega$- limit sets of any admissable, eternal forced mean curvature flow, $[e_t]$, of bounded type and with forcing term $F$ each consist of single points. These points will henceforth be referred to as the {\bf end-points} of $[e_t]$ at plus- and minus- infinity. The results of this section are thus summarised by
\proclaim{Theorem \procref{EndPoints}, {\bf End-points.}}
\noindent If $\Cal{F}$ is of Morse type, and if $[e_t]$ is an admissible, forced mean curvature flow of bounded type with forcing term, $F$, then there exist admissable critical points, $[e_\pm]$, of $\Cal{F}$ such that
$$\eqalign{
\mlim_{t\rightarrow-\infty} [e_t] &= [e_-],\ \text{and}\cr
\mlim_{t\rightarrow+\infty} [e_t] &= [e_+].\cr}
$$
\endproclaim
\newhead{Admissability}[Admissability]
\newsubhead{Pinching}[Pinching]
Consider an Alexandrov embedding, $e:S^d\rightarrow\Bbb{R}^{d+1}$. Denote its shape operator by $A$, and let the subscript ``;'' denote covariant differentiation with respect to the Levi-Civita convariant derivative that it induces over $S^d$. Recall the generalised Simons' formula.
\proclaim{Lemma \nextprocno, {\bf Generalised Simons' formula.}}
\noindent For all $p$, $q$,
$$
A_{pp;qq} = A_{qq;pp} + A^2_{pp}A_{qq} - A_{pp}A^2_{qq}.
$$
\endproclaim
\proclabel{GeneralisedSimonsIdentity}
\proof Indeed, by the Codazzi-Mainardi equations, $A_{ij;k}$ is symmetric under all permutations of the indices. Thus, recalling the definition of curvature,
$$\eqalign{
A_{ij;kl}
&= A_{kj;il}\cr
&= A_{kj;li} - R_{lik}{}^pA_{pj} - R_{lij}{}^pA_{kp}\cr
&= A_{lk;ji} + R_{ilk}{}^pA_{pj} + R_{ilj}{}^pA_{kp},\cr}
$$
where $R$ here denotes the Riemann curvature tensor of the metric induced over $S^d$ by the Alexandrov embedding, $e$. In particular, by Gauss' equation,
$$
R_{ijkl} = A_{il}A_{jk} - A_{ik}A_{jl},
$$
so that
$$
A_{ij;kl} = A_{lk;ji} + A^2_{ij}A_{lk} - A_{lk}A^2_{lj} + A^2_{ik}A_{lj} - A_{ij}A^2_{lk},
$$
and the result now follows upon substituting $i=j=p$ and $k=l=q$.\qed
\medskip
Let $0<\kappa_1\leq...\leq\kappa_d$ be the principal curvatures of $e$. Ideally, the maximum principle would be applied to these functions. However, this would require that they be at least twice differentiable, and since they are, at best, only Lipschitz continuous, they will be smoothed out as follows. Fix a point $(t,x)\in\Bbb{R}\times S^d$, and let $\xi_1,...,\xi_d$ be an orthonormal basis of principal directions of $e_t$ at $x$ corresponding to the principal curvatures $\kappa_1,...,\kappa_d$ respectively. Extend this basis to an orthonormal frame defined in a neighbourhood of $(t,x)$ which is parallel along $S^d$ at $x$, and define the smooth functions $a_1,...,a_d$ in this neighbourhood by $a_k:=\langle A_t \xi_k, \xi_k\rangle$, where $A_t$ is the shape operator of the embedding $e_t$.
\proclaim{Lemma \nextprocno}
\noindent At the point $(t,x)$, for all $k$,
$$
\left(\partial_t - \frac{1}{d}\Delta\right)a_k
=DF(N_t)a_k - a_k^2 F - D^2F(\xi_k,\xi_k) + \frac{1}{d}a_k\opTr(A^2).
$$
\endproclaim
\proclabel{DerivativeOfPrincipalCurvatures}
\proof Let $e:\Bbb{R}\times S^d\rightarrow\Bbb{R}^{d+1}$ be a representative of the flow, $[e_t]$, chosen such that
$$
\partial_t e_t = (H_t - F\circ e_t)N_t.
$$
For fixed vectors, $X$ and $Y$, tangent to $S^d$,
$$\eqalign{
\partial_t (e^*g)(X,Y) &= \partial_t\langle e_*X,e_*Y\rangle\cr
&=\langle\nabla_{\partial_t}e_*X,e_*Y\rangle + \langle e_*X,\nabla_{\partial_t}e_*Y\rangle\cr
&=\langle\nabla_Xe_*\partial_t,e_*Y\rangle + \langle e_*X,\nabla_Ye_*\partial_t\rangle\cr
&=2(H_t - F\circ e_t)\langle A_t X,Y\rangle.\cr}
$$
Since $e^*g(\xi_k,\xi_k)$ is constant, $\partial_t((e^*g)(\xi_k,\xi_k))=0$, and so
$$
e^*g(\partial_t\xi_k,\xi_k) = (H_t - F\circ e_t)\langle A_t\xi_k,\xi_k\rangle,
$$
and combining these relations yields
$$
\partial_t a_k = \langle (\partial_t A)\xi_k,\xi_k\rangle.
$$
However, as is well known (c.f. \cite{Mantegazza}),
$$
\partial_t A = -A_t^2(F\circ e_t-H_t) - \opHess(F\circ e_t - H_t),
$$
where $\opHess$ here denotes the Hessian of functions defined over $S^d$, so that
$$
\partial_t a_k = -a_k^2(F\circ e_t - H_t) - \opHess(F\circ e_t)(\xi_k,\xi_k) + \opHess(H_t)(\xi_k,\xi_k).
$$
Furthermore,
$$
\opHess(F\circ e_t)(\xi_k,\xi_k) = D^2F(\xi_k,\xi_k) - DF(N_t)a_k,
$$
and, by the generalised Simons formula,
$$\eqalign{
\opHess(H_t)(\xi_k,\xi_k) &= \frac{1}{d}\sum_{i=1}^d A_{ii;kk}\cr
&= \frac{1}{d}\sum_{i=1}^d\left(A_{kk;ii} + A_{ii}^2A_{kk} - A_{kk}^2A_{ii}\right)\cr
&= \frac{1}{d}\Delta a_k + \frac{1}{d}a_k\opTr(A_t^2) - a_k^2 H_t.\cr}
$$
The result follows upon combining these relations.\qed
\proclaim{Lemma \nextprocno}
\noindent Let $\lambda\geq 1$, $c>0$ be such that
$$
\frac{(\lambda-1)\lambda}{\lambda+1} > \frac{\|D^2F\|_{L^\infty}}{c^2F_-}.
$$
If $[e_t]$ is $c$-convex for all $t\geq 0$, and if $[e_0]$ is pointwise $\lambda$-pinched, then $[e_t]$ is pointwise strictly $\lambda$-pinched for all $t>0$.
\endproclaim
\proclabel{PreservationOfPinching}
\proof Consider the function $\phi:\Bbb{R}\rightarrow\Bbb{R}$ given by
$$
\phi(t) := \opSup_{x\in S^d}\frac{\kappa_d(t,x)}{\kappa_1(t,x)}.
$$
Since this function is Lipschitz continuous, it is almost everywhere differentiable, its derivative is locally $L^1$, and it is equal to the integral of its derivative. Let $\phi$ be differentiable at the point $t$, and suppose, without loss of generality, that $\phi(t)=\lambda$. Let $x\in S^d$ maximise $\kappa_d/\kappa_1$ and let the functions $a_1$ and $a_d$ be defined near $(t,x)$ as above. In particular, $a_d/a_1$ is smooth, $a_d/a_1\leq\phi$ and $a_d(t,x)/a_1(t,x)=\phi(t)$, so that, at the point $(t,x)$,
$$\eqalign{
\partial_t\opLog(a_d) - \partial_t\opLog(a_1) &= \partial_t\opLog(\phi),\cr
\nabla\opLog(a_d) - \nabla\opLog(a_1) &= 0,\ \text{and}\cr
\Delta\opLog(a_d) - \Delta\opLog(a_1) &\leq 0.\cr}
$$
Thus, by Lemma \procref{DerivativeOfPrincipalCurvatures},
$$\eqalign{
\partial_t\opLog(\phi)
&=\partial_t\opLog(a_d) - \partial_t\opLog(a_1)\cr
&\leq(a_1-a_d)F + \frac{1}{a_1}D^2F(e_1,e_1) - \frac{1}{a_d}D^2F(e_d,e_d)\cr
&\leq a_1(1-\lambda) F + \left(1+\frac{1}{\lambda}\right)\frac{1}{a_1}\|D^2F\|.\cr}
$$
Since $a_1\geq c$, it follows by hypothesis that this is strictly negative, and this completes the proof.\qed
\medskip
This result also yields a straightforward Hopf-type theorem for eternal forced mean curvature flows.
\proclaim{Theorem \nextprocno, {\bf Hopf Theorem for flows.}}
\noindent Let $[e_t]:\Bbb{R}\rightarrow\Cal{E}(\Bbb{R}^{d+1})$ be an eternal forced mean curvature flow of bounded type with constant forcing term $1$. If there exists $\lambda\geq 1$ such that $[e_t]$ is pointwise $\lambda$-pinched and $\lambda$-non-collapsed for all $t$, then, up to translation, $[e_t]$, is the constant flow given by
$$
e_t(x) = x.
$$
\endproclaim
\proclabel{Hopf}
\proof Let $\lambda_0$ be the least value in $[1,\infty[$ for which $[e_t]$ is pointwise $\lambda_0$-pinched. We claim that $\lambda_0=1$. Indeed, otherwise, via a compactness argument based on Theorem \procref{Compactness}, we may suppose that there exists a point $(t,x)\in\Bbb{R}\times S^d$ such that $e_t$ is $\lambda_0$-pinched at the point $x$, which is absurd by Lemma \procref{PreservationOfPinching}. It follows that, for all $t$, $[e_t]$ is totally umbilic at every point, and is therefore a Euclidean hypersphere. Finally, since the flow, $[e_t]$, is of bounded type, and since it exists for all time, the result now follows by Lemmas \procref{DoesNotContainBigSpheres} and \procref{IsNotContainedInSmallSpheres}.\qed
\newsubhead{Non-collapsing}[NonCollapsing]
The following result is an adaptation of the non-collapsing argument, \cite{Andrews}, of Ben Andrews.
\proclaim{Lemma \nextprocno}
\noindent Let $[e_t]$ be an eternal forced mean curvature flow with forcing term, $f$. Let $\lambda\geq 1$, $c>0$ be such that
$$
\frac{(\lambda-1)\lambda}{(\lambda+1)} > \frac{\|D^2F\|_{L^\infty}}{c^2 F_-}.
$$
If $[e_t]$ is $c$-convex and pointwise strictly $\lambda$-pinched for all $t\geq 0$, and if $[e_0]$ is $\lambda$-non-collapsed, then $[e_t]$ is strictly $\lambda$-non-collapsed for all $t>0$.
\endproclaim
\proclabel{PreservationOfNonCollapsing}
\remark An analogous result also holds for spheres {\sl containing} $[e_0]$. This in turn can be used to provide a more direct proof of a-priori upper bounds of the sort obtained in Lemma \procref{BoundedCurvature}. However, since the estimates obtainable by this means are actually weaker, we shall not discuss this further here.
\medskip
\proof For a continuous function, $\Phi:]-\epsilon,\epsilon[\times S^d\rightarrow\Bbb{R}$, consider the function, $Z:\Bbb{R}\times S^d\times S^d\rightarrow\Bbb{R}$, given by
$$
Z(t,x,y) := \Phi(t,x)\|x-y\|^2 + 2\langle y-x,N_t(x)\rangle,
$$
where the points $x$ and $y$ are here identified with their respective images under the embedding $e_t$. Observe that $[e_t]$ contains the Euclidean hypersphere of curvature $\Phi(t,x)$ which is tangent to $[e_t]$ at the point $e(t,x)$ whenever $Z(t,x,y)\geq 0$ for all $y\in S^d$. Suppose therefore that $Z$ attains a minimum of $0$ at the point $y\neq x$, and that $\Phi$ is smooth in a neighbourhood of $(t,x)$. Let $R$ be the reflection orthogonal to $y-x$ which sends $x$ to $y$. Straightforward geometric arguments yield
$$
RN(x) = N(y).
$$
Let $\partial x_1,...,\partial x_d$ be an orthonormal basis of $TS^d$ at $x$ and let $\partial y_1,...,\partial y_d$ be its image under the action of $R$. That is, for all $k$,
$$
R \partial x_k = \partial y_k.
$$
Observe that $\partial y_1,...,\partial y_k$ is also an orthonormal basis of $TS^d$ at $y$. Extend both $(\partial x_k)$ and $(\partial y_k)$ to frames defined in neighbourhoods of $x$ and $y$ respectively which are parallel along all geodesics in $S^d$ leaving these points. Consider now the operator
$$
\Delta^D := \sum_{k=1}^d(\partial x_k + \partial y_k)^2.
$$
The result will follow from the positivity of
$$
\left(\partial_t - \frac{1}{d}\Delta^D\right)Z,
$$
for a suitable choice of the function $\Phi$. Indeed, first, we have
$$\eqalign{
D_{\partial x_i + \partial y_i}(y-x) &= \partial y_i - \partial x_i,\cr
\frac{1}{d}\Delta^D(y-x) &= H_xN_x - H_yN_y,\cr
\partial_t(y-x) &= (F_y-H_y)N_y - (F_x-H_x)N_x,\cr
\left(\partial_t - \frac{1}{d}\Delta^D\right)(y - x) &= F_yN_y - F_xN_x.\cr}
$$
Likewise,
$$\eqalign{
D_{\partial x_i + \partial y_i}N_x &= A_i(\partial x_i)\cr
\frac{1}{d}\Delta^D N_x &= \nabla^\Sigma H - \frac{1}{d}\opTr(A^2)N_x,\cr
\partial_t N_x &= -\nabla^\Sigma(F-H),\cr
\left(\frac{\partial}{\partial t} - \frac{1}{d}\Delta^D\right)(y - x) &= -\nabla^\Sigma F + \frac{1}{d}\opTr(A^2)N_x,\cr}
$$
so that
$$\eqalign{
\left(\partial_t - \frac{1}{d}\Delta^D\right)Z
&=\left[\left(\partial_t - \frac{1}{d}\Delta^D\right)\Phi\right]\|y-x\|^2
-\frac{4}{d}\sum_{i=1}^d(\partial x_i\Phi)\langle y-x,\partial y_i-\partial x_i\rangle\cr
&\qquad + 2\Phi\langle F_y N_y - F_x N_x,y-x\rangle - \frac{2}{d}\Phi\sum_{i=1}^d\left\|\partial y_i-\partial x_i\right\|^2\cr
&\qquad + 2\langle F_y N_y - F_x N_x, N_x\rangle + 2\left\langle y-x,-\nabla^\Sigma F + \frac{1}{d}\opTr(A^2)N_x\right\rangle\cr
&\qquad - \frac{4}{d}\sum_{i=1}^d\langle \partial y_i -\partial x_i, A_x(\partial x_i)\rangle.}
$$
However, since $(x,y)$ minimises $Z$,
$$
(\partial x_i \Phi)\|y-x\|^2 - 2\Phi\langle y-x,\partial x_i\rangle + 2\langle y-x,\partial x_i\rangle\kappa_i = 0,
$$
so that
$$\eqalign{
\langle\partial x_i,y-x\rangle &= \frac{(\partial x_i\Phi)\|y-x\|^2}{2(\Phi - \kappa_i)},\cr
\langle\partial y_i,y-x\rangle &= \frac{-(\partial x_i\Phi)\|y-x\|^2}{2(\Phi - \kappa_i)},\cr
\langle\partial y_i -\partial x_i,y-x\rangle &=\frac{-(\partial x_i\Phi)\|y-x\|^2}{(\Phi-\kappa_i)},\cr
\|\partial y_i \partial x_i\|^2 &= \frac{(\partial x_i\Phi)^2\|y-x\|^2}{(\Phi - \kappa_i)^2},\cr}
$$
and combining these relations yields
$$\eqalign{
\left(\partial_t - \frac{1}{d}\Delta^D\right)Z
&=\|y-x\|^2\bigg[\left(\partial_t-\frac{1}{d}\Delta\right)\Phi + \frac{2}{d}\sum_{i=1}^d\frac{(\partial x_i\Phi)^2}{(\Phi-\kappa_i)}\cr
&\qquad+ \Phi^2 F - \frac{1}{d}\Phi\opTr(A^2) - \Phi\partial_N F\bigg]\cr
&\qquad\qquad+ 2\left(F_y - F_x - \langle y-x,\nabla F\rangle\right)\phantom{\bigg[}.\cr}
$$
In particular, if $\Phi:=\lambda a_1$, where $a_1$ is defined as in Section \subheadref{Pinching}, then, by Lemma \procref{DerivativeOfPrincipalCurvatures},
$$
\left(\partial_t - \frac{1}{d}\Delta^D\right)Z
\geqslant\|y-x\|^2\left(\lambda(\lambda-1)\kappa_1^2 F - (\lambda+1)\|D^2F\|\right)>0,
$$
as desired.\qed %
\medskip
\noindent It follows that admissable eternal forced mean curvature flows of bounded type are contained entirely within $\Cal{E}_{\lambda_F}$.
\proclaim{Theorem \procref{Separation}, {\bf Separation.}}
\noindent Let $[e_t]:\Bbb{R}\rightarrow\Cal{E}(T^{d+1})$ be an eternal forced mean curvature flow of bounded type with forcing term $F$. If $[e_t]$ is pointwise $\Lambda_F$-pinched and $\Lambda_F$-non-collapsed for all $t$, then it is pointwise $\lambda_F$-pinched and $\lambda_F$-non-collapsed for all $t$.
\endproclaim
\proof Let $\lambda_0$ be the least value in $[0,\Lambda_F]$ for which $[e_t]$ is pointwise $\lambda_0$-pinched and $\lambda_0$-non-collapsed, and suppose that $\lambda_0\in]\lambda_F,\Lambda_F[$. Via a compactness argument based on Theorem \procref{Compactness}, we may suppose that there exists a point $(t,x)\in\Bbb{R}\times S^d$ such that either $e_t$ is $\lambda_0$-pinched at $x$ or $e_t$ is $\lambda_0$-non-collapsed at $x$. However, by Lemma \procref{CConvexity}, $[e_t]$ is $c$-convex for all $t$, where $c:=F_-/\lambda_0^2$. Furthermore, since $F$ is subcritical, and since $\lambda_F<\lambda_0<\Lambda_F$, it follows that the hypotheses of Lemmas \procref{PreservationOfPinching} and \procref{PreservationOfNonCollapsing} are both satisfied. This yields a contradiction, and the result follows.\qed
\newhead{Singular perturbation}[SingularPerturbation]
\newsubhead{Asymptotic series}[AsymptoticSeries]
Recall that the main objective of this text is the proof of Theorem \procref{ConcentrationParabolicCase}, which describes the asymptotic behaviour of families of admissable, eternal forced mean curvature flows of bounded type as the forcing term tends to infinity in a certain manner. The first step of this proof involves a presentation of the construction of \cite{SmiEMCFII}, adapted to the case where the ambient manifold is a $(d+1)$-dimensional torus. Throughout this and the following section, extensive use will be made of the straightforward fact that the operation of composition by a given smooth function defines smooth maps between H\"older spaces which send bounded sets to bounded sets.
\par
At this stage, various rather technical definitions will be required. First, given a smooth function, $\phi:S^d\rightarrow]-1,\infty[$, denote
$$
\hat{\phi}(x) := (1+\phi(x))x,
$$
so that $\hat{\phi}$ is simply the natural parametrisation of the graph of $\phi$ over $S^d$. In order to define the functionals of interest to us, it will be useful to introduce the terminology of jet spaces.\footnote*{Consider a finite-dimensional manifold, $M$, and given $k>0$ and a point, $p\in M$, consider the equivalence relation, $\sim_{k,p}$, defined over $C^\infty(M)$ such that $f\sim_{k,p}g$ if and only if they coincide up to order $k$ at the point $p$. The space $J^kE$ of {\bf $k$-jets} of functions over $M$ is defined to be the bundle of all equivalence classes of $\sim_{k,p}$ in $C^\infty(M)$ as $p$ ranges over all points of $M$. Thus, $J^0M$ is simply the trivial bundle $\Bbb{R}\times M$, $J^1M$ is the fibrewise product of $J^0M$ with $T^*M$, and so on.} For all $k$, $J^kS^d$ will denote the bundle of $k$-jets of real-valued functions over $S^d$. The function $H:J^2S^d\rightarrow\Bbb{R}$ will be defined such that, for all $\phi$, and for all $x$, $H(J^2\phi(x))$ will be the mean curvature of the embedding, $\hat{\phi}$, at the point, $x$. Likewise, the function $N:J^1S^d\rightarrow S^d$ will be defined such that, for all $\phi$, and for all $x$, $N(J^1\phi(x))$ will be the outward-pointing, unit, normal vector of the embedding, $\hat{\phi}$, at the point, $x$. $H$ and $N$ will be both smooth (in fact, analytic) functions defined over the finite-dimensional manifolds, $J^2S^d$ and $J^1S^d$, respectively. Furthermore, the Taylor series of $H$ about $0$ will be (c.f. \cite{SmiEMCFII} and \cite{Ye})
$$
H(J^2\phi) = 1 - \frac{1}{d}(d+\Delta)\phi + \sum_{m=2}^\infty P_m(J^2_x\phi),\eqnum{\nexteqnno[TaylorSeriesOfH]}
$$
where, for all $m$, $P_m$ is a homogeneous polynomial of order $m$. Likewise, the Taylor series of $N$ about $0$ will be (c.f. \cite{SmiEMCFII} and \cite{Ye})
$$
N(J^1\phi) = x\left(1 + \sum_{m=2}^\infty P_{1,m}(J_x^1\phi)\right) + \left(-\nabla\phi + \sum_{m=2}^\infty P_{2,m}(J_x^1\phi)\right),\eqnum{\nexteqnno[TaylorSeriesOfN]}
$$
where, for all $m$, $P_{1,m}$ and $P_{2,m}$ are homogeneous polynomials of order $m$.
\par
Now consider a smooth function, $f:T^{d+1}\rightarrow\Bbb{R}$, of Morse-Smale type, and let $\gamma:\Bbb{R}\rightarrow T^{d+1}$ be one of its complete gradient flows (c.f. Section \subheadref{MorseHomology}). For $\kappa>0$ and for smooth functions, $X:\Bbb{R}\rightarrow\Bbb{R}^{d+1}$ and $\phi:\Bbb{R}\times S^d\rightarrow]-1/\kappa,\infty[$, the smooth function, $\Phi(\kappa,X,\phi):\Bbb{R}\times S^d\rightarrow T^{d+1}$, will be defined such that
$$
\Phi(\kappa,X,\phi)(t,x) := \gamma(t) + \kappa X(t) + \kappa(1 + \kappa^2\phi(t,x))x,\eqnum{\nexteqnno[FlowFunction]}
$$
so that, heuristically, $\Phi(\kappa,X,\phi)$ will be a smooth family of embedded spheres each of radius approximately $\kappa$, with centres moving along $\gamma$, though displaced slightly by the vector field, $\kappa X$.
\par
For all $\kappa>0$, the forcing term, $F_\kappa:T^{d+1}\rightarrow\Bbb{R}$, will be defined by
$$
F_\kappa := \frac{1}{\kappa}(1 - \kappa^2 f).
$$
We now study the conditions to be satisfied by the triplet $(\kappa,X,\phi)$ in order for the family $\Phi(\kappa,X,\phi)$ to define an eternal forced mean curvature flow with forcing term, $F_\kappa$. We require a further extension of our terminology of jet spaces. Thus, for all $k$, $J^k(\Bbb{R},\Bbb{R}^{d+1})$ will denote the bundle of $k$-jets of smooth $\Bbb{R}^{d+1}$-valued functions over $\Bbb{R}$. For all $k$, and for any smooth function, $\psi:\Bbb{R}\times S^d\rightarrow\Bbb{R}$, $J^k_\opin\psi$ will denote the inhomogeneous jet consisting of all derivatives of $\psi$ of order $i$ in $x$ and order $j$ in $t$ for all $i+2j\leq 2k$. For all $k$, $J^k_\opin(\Bbb{R}\times S^d)$ will denote the bundle of inhomogeneous $k$-jets of real-valued functions over $\Bbb{R}\times S^d$. Observe that $J^k_\opin(\Bbb{R}\times S^d)$ is also a bundle over $\Bbb{R}$, and, for all $k$, $J^k$ will denote its fibrewise product over $\Bbb{R}$ with $J^k(\Bbb{R},\Bbb{R}^{d+1})$. In other words, a typical element of $J^k$ will have the form $(J^kX(t),J^k_\opin\phi(t,x))$, for some smooth functions, $X:\Bbb{R}\rightarrow\Bbb{R}^{d+1}$ and $\phi:\Bbb{R}\times S^d\rightarrow\Bbb{R}$, and for some point $(t,x)\in\Bbb{R}\times S^d$.
\par
The function $\Psi:]0,\infty[\times J^1\rightarrow\Bbb{R}$ will be defined such that, for all $\kappa\in]0,\infty[$, for all smooth functions, $X:\Bbb{R}\rightarrow\Bbb{R}^{d+1}$ and $\phi:\Bbb{R}\times S^d\rightarrow\Bbb{R}$, and for any point $(t,x)\in\Bbb{R}\times S^d$,
$$\eqalign{
\Psi(\kappa,J^1X(t),J^1_\opin \phi(t,x))
&:=\kappa^2\langle D_t\Phi(\kappa,X,\phi)(t,x),N(\kappa^2J^1_x\phi(t,x))\rangle+\cr
&\qquad\qquad\qquad\frac{1}{\kappa}H(\kappa^2J_x^2\phi(t,x)) - F_\kappa(\Phi(\kappa,X,\phi)(t,x)),}\eqnum{\nexteqnno[DefinitionOfPsi]}
$$
where $J_x^1\phi$ and $J_x^2\phi$ here denote respectively the first and second order jets of derivatives of $\phi$ in the $x$ direction. As with $H$ and $N$, $\Psi$ will be a smooth function defined over the finite-dimensional manifold, $]0,\infty[\times J^1$. By \eqnref{ForcedMCFEquation}, we see that this function has been constructed precisely in order that, given $\kappa>0$, and smooth functions, $X:\Bbb{R}\rightarrow\Bbb{R}^{d+1}$ and $\phi:\Bbb{R}\times S^d\rightarrow\Bbb{R}$, the family $\Phi(\kappa,X,\phi)$ is an eternal forced mean curvature flow with forcing term, $F_\kappa$, if and only if $\Psi(\kappa,J^1X(t),J^1_\opin\phi(t,x))$ vanishes for all $(t,x)\in\Bbb{R}\times S^d$.
\par
Solutions of the forced mean curvature flow equation will now be obtained in the standard manner by first constructing formal solutions which are then perturbed to exact solutions via the inverse function theorem. To this end, the following terminology of asymptotic series will also prove useful. Consider an arbitrary, finite-dimensional vector bundle, $E$, over some base, $B$. Given a smooth function, $\phi:]0,\infty[\times E\rightarrow\Bbb{R}$, and a sequence, $(\phi_i)$, of smooth functions where, for all $i$, $\phi_i:E^{\otimes i}\rightarrow\Bbb{R}$, and a formal power series, $\xi(\kappa,x)\sim\sum_{i=0}^\infty \kappa^i\xi_i(x)$, taking values in $E$, the expression
$$
\phi(\kappa,\xi)\sim\sum_{i=0}^\infty \kappa^i\phi_i(\xi_{0,x},...,\xi_{i,x})\eqnum{\nexteqnno[AsymptoticSeries]}
$$
will mean that, for all $m\geq 0$, there exists a smooth function, $R_m:[0,\infty[\times E^{\otimes m}\rightarrow\Bbb{R}$, such that
$$
\phi\left(\kappa,\sum_{i=0}^m \kappa^i\xi_i(x)\right) = \sum_{i=0}^m\kappa^i\phi_i\left(\xi_0(x),...,\xi_i(x)\right) + \kappa^{m+1}R_m(\kappa,\xi_0(x),...,\xi_m(x)).
$$
The series $\sum_{i=0}^\infty\kappa^i\phi_i$ will be called the {\bf asymptotic series} of the function, $\phi$. It will be of fundamental importance here that the remainder term, $R_m$, be smooth at the point $t=0$. Indeed, otherwise, this relation would be trivially satisfied for any pair of sequences, $(\phi_i)$ and $(\xi_i)$, making it of little use to us.
\proclaim{Lemma \nextprocno}
\noindent There exists a sequence, $(\Psi_i)$, of functions, which are polynomials in all but the first variable, such that, for all formal series, $X\sim\sum_{i=0}^\infty\kappa^i X_i$ and $\phi\sim\sum_{i=0}^\infty\kappa^i\phi_i$,
$$\multiline{
\Psi(\kappa,X,\phi)\sim&
-\frac{1}{2}\kappa^3\opHess(f)(\gamma(t))(X_0,X_0) - \frac{1}{2}\kappa^3\opHess(f)(\gamma(t))(x,x)\cr
& +\kappa \sum_{i=0}^\infty\kappa^i\bigg[\langle D_t X_{i-2} - \opHess(f)(\gamma(t))X_{i-2},x\rangle\cr
&\qquad + \left(\kappa^4 D_t - \frac{1}{d}(d+\Delta)\right)\phi_i\cr
&\qquad + \langle\nabla f(\gamma(t)),X_{i-1}\rangle + \opHess(f)(\gamma(t))(X_{i-2},X_0)\vphantom{\bigg)}\cr
&\qquad + \Psi_m(\gamma,D_t\gamma,x,J^1_tX_0,...,J^1_tX_{i-3},\vphantom{\bigg)}\cr
&\qquad\qquad \kappa^4 D_t\phi_0,...,\kappa^4 D_t\phi_{i-2},J_x^2\phi_0,...,J_x^2\phi_{i-2})\bigg],}\eqnum{\nexteqnno[AsymptoticSeriesOfPsi]}
$$
where, by convention, $X_i$ and $\phi_i$ are taken to be equal to $0$ for $i$ negative. Furthermore,
$$\eqalign{
\Psi_0 &= f(\gamma(t)),\ \text{and}\cr
\Psi_1 &= 0.\cr}
$$
\endproclaim
\proof By \eqnref{TaylorSeriesOfH}, there exists a sequence, $(H_i)$, of polynomials such that
$$
H(\kappa^2J_x^2\phi) \sim 1-\kappa^2\sum_{i=0}^\infty\frac{1}{d}(d + \Delta)\kappa^i\phi_i + \kappa^4\sum_{i=0}^\infty\kappa^i H_i(J_x^2\phi_0,...,J_x^2\phi_i).
$$
Likewise, by \eqnref{TaylorSeriesOfN}, there exist sequences, $(N_{1,i})$ and $(N_{2,i})$, of polynomials such that
$$
N(x,\kappa^2 J_x^1\psi) = x(1 + N_1(J_x^2\psi)) + \kappa^2 N_2(J_x^1\psi),
$$
where
$$\eqalign{
N_1(\kappa^2 J_x^1\psi) &\sim \kappa^4\sum_{i=0}^\infty\kappa^i N_{1,i}(J_x^1\psi_0,...,J_x^1\psi_i),\ \text{and}\cr
N_2(\kappa^2 J_x^1\psi) &\sim -\sum_{i=0}^\infty\kappa^i\nabla\psi_i + \kappa^2\sum_{i=0}^\infty\kappa^i N_{2,i}(J_x^1\psi_0,...,J_x^1\psi_i).}
$$
Finally, by Taylor's theorem, there exists a sequence, $(F_i)$, of functions, which are polynomials in all but the first variable, such that
$$\eqalign{
F_\kappa(\Phi(\kappa,X,\phi))
&\sim \frac{1}{\kappa} - \kappa f(\gamma(t)) - \kappa^2\langle\nabla f(\gamma(t)),x\rangle\cr
&\qquad -\kappa^2\sum_{i=0}^\infty\kappa^i\langle\nabla f(\gamma(t)),X_i\rangle\cr
&\qquad -\frac{1}{2}\kappa^3\langle\opHess(f)(\gamma(t))x,x\rangle
-\frac{1}{2}\kappa^3\sum_{i,j=0}^\infty\kappa^{i+j}\langle\opHess(f)(\gamma(t))X_i,X_j\rangle\cr
&\qquad -\kappa^3\sum_{i=0}^\infty\kappa^i\langle\opHess(f)(\gamma(t))X_i,x\rangle\cr
&\qquad +\kappa^4\sum_{i=0}^\infty F_i(\gamma(t),x,J_t^0X_0,...,J_t^0X_i,J_x^0\phi_0,...,J_x^0\phi_i).\cr}
$$
The result now follows upon combining these relations with the definitions \eqnref{FlowFunction} and \eqnref{DefinitionOfPsi} of $\Phi$ and $\Psi$ respectively.\qed
\newsubhead{Spherical harmonics and formal solutions}[SphericalHarmonicsAndFormalSolutions]
Consider the standard Laplace-Beltrami operator, $\Delta$, of functions over $S^d$. For all integer $m\geq 0$, the vector space, $\Cal{H}_m:=\Cal{H}_m(S^d)$, of eigenfunctions of $\Delta$ with eigenvalue $m(m+d-1)$ will be referred to as the space of $m$'th order, d-dimensional {\bf spherical harmonics}. Since every eigenfunction of $\Delta$ is a spherical harmonic of some order, the spherical harmonics together form an orthonormal basis of $L^2(S^d)$. For all $m$, every $m$'th order spherical harmonic is the restriction to $S^d$ of a unique homogeneous, harmonic polynomial of degree $m$. Conversely, the restriction to $S^d$ of any polynomial of degree $m$ decomposes uniquely as a sum of spherical harmonics of orders at most $m$. In addition, if this polynomial is even (resp. odd), then every non-trivial component of this decomposition has even (resp. odd) order.
\par
The mean curvature flow equation, $\Psi(\kappa,J^1X(t),J^1_\opin\phi(t,x))=0$, introduced in the previous section, exhibits two distinct asymptotic behaviours as the parameter, $\kappa$, tends to $0$. This will be resolved by rescaling the space, $\Cal{H}_1$, of first-order spherical harmonics, as the parameter, $\kappa$, tends to $0$ (c.f. Section \subheadref{ExactSolutions}, below). To this end, we will denote
$$
\hat{\Cal{H}} := \oplus_{m\neq 1}\Cal{H}_m,
$$
and $\Pi:L^2(S^d)\rightarrow\Cal{H}_1$ and $\Pi^\perp:L^2(S^d)\rightarrow\hat{\Cal{H}}$ will denote the orthogonal projections. In like manner, for all $k$, we will denote
$$
\hat{\Cal{H}}_k := \oplus_{m\neq 1,m\leq k}\Cal{H}_m.
$$
In addition, $C^\infty_\opbdd(\Bbb{R}\times S^d)$ will denote the space of smooth functions, $f:\Bbb{R}\times S^d\rightarrow\Bbb{R}$, whose derivatives to all orders are bounded, and $\hat{C}^\infty_\opbdd(\Bbb{R}\times S^d)$ will denote the subspace consisting of those functions, $f:\Bbb{R}\times S^d\rightarrow\Bbb{R}$, such that, for all $t$,
$$
\Pi(f(t,\cdot)) = 0.
$$
Similarly, for any finite-dimensional vector space, $E$, $C^\infty_\opbdd(\Bbb{R},E)$ will denote the space of smooth functions, $f:\Bbb{R}\rightarrow E$, whose derivatives to all orders are bounded. In particular, for all $k\neq 1$, $C^\infty_\opbdd(\Bbb{R},\hat{\Cal{H}}_k)$ canonically identifies with a subspace of $\hat{C}^\infty_\opbdd(\Bbb{R}\times S^d)$.
\par
Consider now the operators that arise in the asymptotic series, \eqnref{AsymptoticSeriesOfPsi}, of $\Psi$. First, the operator, $L:C^\infty_\opbdd(\Bbb{R},\Bbb{R}^{d+1})\rightarrow C^\infty_\opbdd(\Bbb{R},\Bbb{R}^{d+1})$, is defined by
$$
L:=D_t+\opHess(f)(\gamma(t)).\eqnum{\nexteqnno[DefinitionOfL]}
$$
Since $f$ is of Morse type, this operator is of Fredholm type with Fredholm index equal to the difference of the Morse indices of its end-points (c.f. \cite{Schwarz}). Furthermore, elements of its kernel decay exponentially at $\pm\infty$, so that the $L^2$-orthogonal complement of this kernel in $C^\infty_\opbdd(\Bbb{R},\Bbb{R}^{d+1})$ is well-defined. Finally, since $f$ is of Morse-Smale type, this operator is surjective, and therefore has a unique left-inverse, $K$, which sends $C^\infty_\opbdd(\Bbb{R},\Bbb{R}^{d+1})$ into $\opKer(L)^\perp$.
\par
Next, for all $\kappa>0$, the operator, $P_\kappa:C^\infty_\opbdd(\Bbb{R}\times S^d)\rightarrow C^\infty_\opbdd(\Bbb{R}\times S^d)$, is defined by
$$
P_\kappa:=\kappa^4 D_t -\frac{1}{d}\left(d + \Delta\right).\eqnum{\nexteqnno[DefinitionOfP]}
$$
It follows from the classical theory of parabolic operators that $P_\kappa$ defines an invertible map from $\hat{C}^\infty_\opbdd(\Bbb{R}\times S^d)$ to itself. Let $Q_\kappa$ denote its inverse. Observe, furthermore, that, for all $k\neq 1$, both $P_\kappa$ and $Q_\kappa$ restrict to linear isomorphisms of $C^\infty_\opbdd(\Bbb{R},\hat{\Cal{H}}_k)$.
\par
Formal solutions of the forced mean curvature flow equation are now obtained in a straightforward manner. Oddly, the coefficients of the asymptotic series of $X$ and $\phi$ actually depend on $\kappa$, although, as we will show presently, they are nonetheless uniformly bounded in every norm.
\proclaim{Lemma \nextprocno}
\noindent There exist unique formal series, $X\sim\sum_{i=0}^\infty\kappa^iX_{i,\kappa}$ in $\opKer(L)^\perp$, and $\phi\sim\sum_{i=0}^\infty\kappa^i\phi_{i,\kappa}$ in $\hat{C}^\infty_\opbdd(\Bbb{R}\times S^d)$, formally solving
$$
\Phi(\kappa,J^1X(x),J^1_\opin\phi(t,x)) \sim 0.
$$
Furthermore, for all $i$, there exists $k_i$ such that $\phi_i$ is an element of $C^\infty_\opbdd(\Bbb{R},\hat{\Cal{H}}_{k_i})$.
\endproclaim
\proclabel{FormalSolutions}
\proof Suppose, by induction, that $X_{0,\kappa},...,X_{i,\kappa}$ and $\phi_{0,\kappa},...,\phi_{i,\kappa}$ have been constructed. In particular, the orthogonal projections onto $\hat{\Cal{H}}$ of terms in the series \eqnref{AsymptoticSeriesOfPsi} of order up to and including $(i+1)$ in $\kappa$ all vanish. Likewise, the orthogonal projections onto $\Cal{H}_1$ of terms in this series of order up to and including $(i+3)$ in $\kappa$ also all vanish. Define
$$\eqalign{
\phi_{i+1,\kappa} &:= -Q_\kappa \left(\langle\nabla f(p),X_{i,\kappa}\rangle + \opHess(f)(p)(X_{i-1,\kappa},X_{0,\kappa})\right.\cr
&\qquad\qquad + \Pi_1^\perp\Psi_i(D_t\gamma,x,J_t^1X_{0,\kappa},...,J_t^1X_{i-2,\kappa},\cr
&\qquad\qquad\qquad\qquad \left.\kappa^4D_t\phi_{0,\kappa},...,\kappa^4D_t\phi_{i-1,\kappa},J_x^2\phi_{0,\kappa},...,J_x^2\phi_{i-1,\kappa})\right).\cr}
\eqnum{\nexteqnno[InductiveDefinitionOfPhi]}
$$
In particular, since $\Psi_i$ is a polynomial in all but the first variable, it follows by the inductive hypothesis that $\phi_{i+1,\kappa}$ is an element of $C^\infty_\opbdd(\Bbb{R},\hat{\Cal{H}}_{k_{i+1}})$, for some $k_{i+1}$ independent of $\kappa$. Define
$$\eqalign{
X_{i+1,\kappa} &:= -L\left(\Pi_1\Psi_i(D_r\gamma,x,J_t^1X_{0,\kappa},...,J_t^1X_{i,\kappa},\right.\cr
&\qquad\qquad\qquad \left.\kappa^4 D_t\phi_{0,\kappa},...,\kappa^4 D_t\phi_{i+1,\kappa},J_x^2\phi_{0,\kappa},...,J_x^2\phi_{i+1,\kappa})\right).\cr}
\eqnum{\nexteqnno[InductiveDefinitionOfX]}
$$
With $\phi_{i+1,\kappa}$ defined as above, the orthogonal projection onto $\hat{\Cal{H}}$ of the term in the series \eqnref{AsymptoticSeriesOfPsi} of order $(i+2)$ in $\kappa$ vanishes. Likewise, with $X_{i+1,\kappa}$ defined as above, the orthogonal projection onto $\Cal{H}_1$ of the term in this series of order $(i+4)$ in $\kappa$ also vanishes. Furthermore, by definition of $Q_\kappa$ and $L$, $\phi_{i+1,\kappa}$ and $X_{i+1,\kappa}$ are the only functions with these properties, and the result now follows.\qed
\newsubhead{The $\kappa$-dependent rescaled time variable}[TheKappaDependentRescaledTimeVariable]
Consider a finite-dimensional vector space, $E$, and a linear operator, $M:C^1_\opbdd(\Bbb{R},E)\rightarrow C^0_\opbdd(\Bbb{R},E)$, of the form
$$
M f  = D_t f + A(t)f,
$$
where $A:\Bbb{R}\rightarrow\opEnd(E)$ is a smooth function whose derivatives to all orders are bounded, and such that
$$
\mlim_{t\rightarrow\pm\infty}A(t) = A_\pm,
$$
for some invertible operators $A_\pm$. These hypotheses imply that $M$ is of Fredholm type (c.f. \cite{Schwarz}). Furthermore, all derivatives of all elements of its kernel decay exponentially at plus- and minus- infinity, and, in particular, the $L^2$-orthogonal projection, $\Pi$, from $C^0_\opbdd(\Bbb{R},E)$ onto this kernel is well-defined. Let $\|\cdot\|_0$ denote the uniform norm over $C^0_\opbdd(\Bbb{R},E)$. Straightforward integral estimates yield
\proclaim{Lemma \nextprocno}
\noindent There exists a constant, $B>0$, such that, for all $f\in C^1_\opbdd(\Bbb{R},E)$,
$$
\|f\|_0, \|D_tf\|_0 \leq B\left(\|Mf\|_0 + \|\Pi f\|_0\right).
$$
\endproclaim
As can be seen from \eqnref{DefinitionOfL} and \eqnref{DefinitionOfP}, different components of the flow, $\Phi(\kappa,X,\phi)$, evolve at different rates, depending on the parameter, $\kappa$. For this reason, it will be useful to introduce the {\bf $\kappa$-dependent, rescaled time variable}, $s:=\kappa^{-4}t$, so that the operator of differentiation with respect to $s$ will then be given by $D_s:=\kappa^4D_t$.
\par
For all $\alpha\in]0,1]$, the $\kappa$-dependent {\bf H\"older difference operator} of order $\alpha$ will be given by
$$
\delta_s^\alpha f(t) := \msup_{\left|h\right|\leq 1}\frac{\left|f(t+\kappa^4 h) - f(t)\right|}{h^\alpha}.
$$
Likewise, for all $(k,\alpha)$, the $\kappa$-dependent {\bf H\"older norm} of order $(k,\alpha)$ will be given by
$$
\|f\|_{k,\alpha,\kappa} :=  \sum_{i=0}^k\|D_s^if\|_0 + \|\delta_s^\alpha D_s^k f\|_0.
$$
Finally, the {\bf H\"older space} of order $(k,\alpha)$ will be defined by
$$
C^{k,\alpha}(\Bbb{R},E) := \left\{ f\in C^k(\Bbb{R},E)\ |\ \|f\|_{k,\alpha,\kappa} < \infty\right\}.
$$
Observe that, for any given $(k,\alpha)$, the $\kappa$-dependent H\"older norms of order $(k,\alpha)$ are all pairwise uniformly equivalent, so that the space, $C^{k,\alpha}(\Bbb{R},E)$, defined as above, is indeed independent of $\kappa$.
\par
An important fact to be used in the sequel is that the norm of the Green's operator of $M$ is actually uniformly bounded, independent of the parameter, $\kappa$. Indeed,
\proclaim{Lemma \nextprocno}
\noindent For all $(k,\alpha)$, there exists $B>0$ such that, for all $\kappa\in]0,1]$, and for all $f\in C^{k+1,\alpha}(\Bbb{R},E)$,
$$
\|f\|_{k,\alpha,\kappa},\|D_tf\|_{k,\alpha,\kappa} \leq B\left(\|Mf\|_{k,\alpha,\kappa} + \|\Pi f\|_0\right).
$$
\endproclaim
\proclabel{UniformNormOfK}
\proof Observe that, for all $(k,\alpha)$, there exists $B_{k,\alpha}>0$ such that, for all $\kappa\in]0,1]$,
$$
\|\delta_s^\alpha D_s^k A\|_0 \leq B_{k,\alpha}\kappa^{4(k+\alpha)} \leq B_{k,\alpha}.\eqnum{\nexteqnno[BoundsForA]}
$$
Let $h_1,...,h_m$ be an $L^2$-orthonormal basis of $\opKer(M)$. Recall (c.f. \cite{Schwarz}), that, for all $i$, and for all $k$, $D_t^kh_i$ decays exponentially at plus- and minus- infinity, so that, upon increasing $B_{k,\alpha}$ if necessary, we may suppose that,
$$
\int_\Bbb{R}\|\delta_s^\alpha D_s^k h_i\| dt \leq B_{k,\alpha}\kappa^{4(k+\alpha)} \leq B_{k,\alpha}.\eqnum{\nexteqnno[BoundsForH]}
$$
Now suppose, by induction, that the result holds for all $0\leq k\leq l$. Then,
$$
\|D_s^{l+1}f\|_0 \leq \|MD_s f\|_{l,\alpha,\kappa} + \|\Pi D_s f\|_0.
$$
However,
$$
M D_s f = D_s Mf + (D_sA)f,
$$
so that, by \eqnref{BoundsForA}, and the inductive hypothesis,
$$\eqalign{
\|MD_s f\|_{l,\alpha,\kappa}
&\leq B_1\left(\|D_sM f\|_{l,\alpha,\kappa} + \|f\|_{l,\alpha,\kappa}\right)\cr
&\leq B_2\left(\|Mf\|_{l+1,\alpha,\kappa} + \|\Pi(f)\|_0\right),\cr}
$$
for some constants, $B_1,B_2>0$. Likewise, by \eqnref{BoundsForH}, and the inductive hypothesis again, for all $i$,
$$\eqalign{
\left|\int_\Bbb{R}(D_s f) h_i dt\right|
&=\left|\int_\Bbb{R} f(D_s h_i)dt\right|\cr
&\leq B_3\|f\|_0\cr
&\leq B_4\left(\|Mf\|_{l+1,\alpha,\kappa} + \|\Pi(f)\|_0\right),\cr}
$$
for some constants, $B_3,B_4>0$. Upon combining these relations, it follows that
$$
\|D_s^{l+1}f\|_0 \leq B_5\left(\|Mf\|_{l+1,\alpha,\kappa} + \|\Pi f\|_0\right),
$$
for some constant, $B_5>0$. Furthermore, using \eqnref{BoundsForA} and the inductive hypothesis once again,
$$\eqalign{
\|D_tD_s^{l+1}f\|_0 &= \|D_s^{l+1}D_tf\|_0\cr
&\leq \|D_s^{l+1}Mf\|_0 + \|D_s^{l+1}Af\|_0\cr
&\leq B_6\left(\|Mf\|_{l+1,\alpha,\kappa} + \|\Pi f\|_0\right),\cr}
$$
for some constant, $B_6>0$. The estimates for $\delta_s^\alpha D_s^{l+1} f$ and $\delta_s^\alpha D_s^{l+1}D_t f$ are obtained in a similar manner, and the result now follows by induction.\qed
\medskip
It will also be useful to introduce $\kappa$-dependent H\"older norms for spaces of functions defined over $\Bbb{R}\times S^d$. Thus, given $\alpha\in]0,1]$, the $\kappa$-dependent {\bf H\"older difference operators} of order $\alpha$ are given by
$$\eqalign{
\delta_x^\alpha f(t,x)
&:= \msup_{y\neq x}\frac{\left|f(t,y) - f(t,x)\right|}{\|x-y\|^\alpha},\ \text{and}\cr
\delta_s^\alpha f(t,x)
&:= \msup_{\left|h\right|\leq 1}\frac{\left|f(t + \kappa^4 h,x) - f(t,x)\right|}{h^\alpha}.\cr}
$$
For all $k\in\Bbb{N}$, $C^k_\opin(\Bbb{R}\times S^d)$ will denote the space of all functions, $f:\Bbb{R}\times S^d\rightarrow\Bbb{R}$ which are continuously differentiable $i$ times in the $x$ direction and $j$ times in the $t$ direction for all $i+2j\leq 2k$. For all $k\in\Bbb{N}$ and for all $\alpha\in ]0,1/2]$, the $\kappa$-dependent {\bf inhomogeneous H\"older norm} of order $(k,\alpha)$ over $C^k_\opin(\Bbb{R}\times S^m)$ will be defined by
$$
\|f\|_{k,\alpha,\kappa,\opin} := \sum_{i+2j\leq 2k}\|D_x^i D_s^j f\|_0
+ \sum_{i+2j=2k}\|\delta_x^{2\alpha}D_x^iD_s^jf\|_0
+ \sum_{i+2j=2k}\|\delta_s^\alpha D_x^iD_s^jf\|_0.
$$
For all $(k,\alpha)$, the {\bf inhomogeneous H\"older space} of order $(k,\alpha)$ will be defined by
$$
C^{k,\alpha}_{\opin}(\Bbb{R}\times S^d) := \left\{ f\in C^k_\opin(\Bbb{R}\times S^d)\ |\ \|f\|_{k,\alpha,\kappa,\opin}<\infty\right\}.
$$
As before, for any given $(k,\alpha)$, the $\kappa$-dependent inhomogeneous H\"older norms of order $(k,\alpha)$ are all pairwise uniformly equivalent, so that the space, $C^{k,\alpha}_{\opin}(\Bbb{R}\times S^d)$, as defined above, is indeed independent of $\kappa$.
\par
Uniform estimates for the formal solutions constructed in the preceeding section now follow.
\proclaim{Lemma \nextprocno}
\noindent For all $(k,\alpha)$, and for all $i$, there exists $B_{k,\alpha,i}>0$ such that, for all $\kappa\in]0,1]$,
$$\eqalign{
\|X_{i,\kappa}\|_{k,\alpha,\kappa} &\leq B_{k,\alpha,i},\ \text{and}\cr
\|\phi_{i,\kappa}\|_{k,\alpha,\kappa,\opin} &\leq B_{k,\alpha,i}.\cr}
$$
\endproclaim
\proclabel{UniformBoundsOnCoefficients}
\proof Indeed, suppose, by induction, that the result holds for all $i\leq j$. Recall that composition by a given smooth function defines smooth maps between H\"older spaces which send bounded sets to bounded sets. Thus, since $\Psi_j$ is itself a smooth function defined over some finite-dimensional manifold, it follows by \eqnref{InductiveDefinitionOfPhi}, \eqnref{InductiveDefinitionOfX} and the inductive hypothesis that
$$\eqalign{
\|P_\kappa\phi_{j+1,\kappa}\|_{k,\alpha,\kappa,\opin}&\leq B_{k,\alpha,j+1},\ \text{and}\cr
\|LX_{j+1,\kappa}\|_{k,\alpha,\kappa}&\leq B_{k,\alpha,j+1},\cr}
$$
for some constant, $B_{k,\alpha,j+1}>0$. The estimate for $X_{j+1,\kappa}$ now follows by Lemma \procref{UniformNormOfK}. However, for all $\kappa$, $P_\kappa = D_s - \frac{1}{d}(d+\Delta)$, so that the norm of the operator, $Q_\kappa$, with respect to the $\kappa$-dependent H\"older norms, $\|\cdot\|_{k,\alpha,\kappa,\opin}$ and $\|\cdot\|_{k+1,\alpha,\kappa,\opin}$, is actually independent of $\kappa$. The estimate for $\phi_{j+1,\kappa}$ therefore follows, and this completes the proof.\qed
\newsubhead{Exact solutions}[ExactSolutions]
Consider now the map $\hat{\Psi}_\kappa:C^{k+1,\alpha}(\Bbb{R},\Bbb{R}^{d+1})\times\hat{C}^{k+1,\alpha}_{\opin}(\Bbb{R}\times S^d)\rightarrow C^{k,\alpha}_{\opin}(\Bbb{R}\times S^d)$ given by
$$
\hat{\Psi}_\kappa(X,\phi)(t,x) := \frac{1}{\kappa}\Psi(\kappa,J^1X(t),J^1_\opin\phi(t,x)).
$$
Since $\hat{\Psi}_\kappa$ is constructed via a combination of differentiation and composition by smooth functions, it defines a smooth map between H\"older spaces. Recall that, as outlined in Section \subheadref{SphericalHarmonicsAndFormalSolutions}, for all $(k,\alpha)$, the space $C^{k,\alpha}_\opin(\Bbb{R}\times S^d)$ naturally decomposes as the direct sum
$$
C^{k,\alpha}_\opin(\Bbb{R}\times S^d) = C^{k,\alpha}(\Bbb{R},\Bbb{R}^{d+1})\oplus\hat{C}^{k,\alpha}_\opin(\Bbb{R}\times S^d).
$$
With respect to this decomposition, the linear isomorphism $T_\kappa:C^{k,\alpha}_\opin(\Bbb{R}\times S^d)\rightarrow C^{k,\alpha}_\opin(\Bbb{R}\times S^d)$ is defined by
$$
T_\kappa := \pmatrix \kappa^{-2}\opId \hfill&\cr &\opId\hfill\cr\endpmatrix.
$$
Using this operator, uniform a-priori bounds for the inverse of the derivative of $\hat{\Psi}$ are obtained.
\proclaim{Lemma \nextprocno}
\noindent For all $(k,\alpha)$, and for all $R>0$, there exists $B>0$ such that, for sufficiently small $\kappa$, and for all $\|X\|_{k+1,\alpha,\kappa} + \|\phi\|_{k+1,\alpha,\kappa,\opin}<R$, the operator, $T_\kappa D\hat{\Psi}_\kappa(X,\phi)$, defines a linear isomorphism from $C^{k+1,\alpha}(\Bbb{R},\Bbb{R}^{d+1})\times\hat{C}^{k+1,\alpha}_\opin(\Bbb{R}\times S^d)$ into $C^{k,\alpha}_\opin(\Bbb{R}\times S^d)$, such that
$$
\left\|\left[T_\kappa D\hat{\Psi}_\kappa(X,\phi)\right]^{-1}\right\| \leq B.
$$
\endproclaim
\proclabel{BoundedInverse}
\proof Indeed, consider first the separate components of the derivative of $\hat{\Psi}_\kappa$. By \eqnref{TaylorSeriesOfH}, the derivative of $\kappa^{-2}H(\kappa^2 J^2_x\phi)$ is given by
$$
\kappa^{-2}DH(\kappa^2J_x^2\phi)\psi = -\frac{1}{d}(d+\Delta)\psi + \kappa^2N_1J_x^2\psi,
$$
where $N_1$ is a multiplication operator of norm bounded in terms of $R$. By \eqnref{TaylorSeriesOfN}, the derivative of $N(\kappa^2 J_x^1\phi)$ is given by
$$
DN(\kappa^2J_x^1\phi)\psi = \kappa^2N_2J_x^1\psi,
$$
where $N_2$ is a multiplication operator of norm bounded in terms of $R$. Thus
$$\eqalign{
\kappa D\langle D_t\Phi_\kappa(X,\phi),N(\kappa^2 J_x^1\phi)\rangle(Y,\psi)
&=\kappa^2\langle D_tY,x\rangle + \kappa^4 D^t\psi\cr
&\qquad + \kappa^2N_3(D_s\psi) + \kappa^4N_4J_x^1\psi + \kappa^4N_5D_tY,\cr}
$$
where $N_3$, $N_4$ and $N_5$ are multiplication operators of norms bounded in terms of $R$. Finally, the derivative of $\frac{1}{\kappa}F(\Phi_\kappa(X,\phi))$ is given by
$$\eqalign{
\frac{1}{\kappa}DF(\Phi_\kappa(X,\phi))
&=-\kappa\langle\nabla f(\gamma(t)), Y\rangle - \kappa^2\langle\opHess(f)(\gamma(t))Y,x\rangle\cr
&\qquad\qquad -\kappa^2\langle\opHess(\gamma(t))Y,X\rangle + \kappa^3N_6Y + \kappa^3 N_7\psi,}
$$
where $N_6$ and $N_7$ are multiplication operators of norms bounded in terms of $R$. Combining these relations, it follows that, with respect to the above decomposition of $C^{k,\alpha}_\opin(\Bbb{R}\times S^d)$, the derivative, $D\hat{\Psi}_\kappa(X,\phi)$, satisfies
$$
T_\kappa D\hat{\Psi}_\kappa(X,\phi) =
\pmatrix L\hfill& M_1\hfill\cr
&P_\kappa\hfill\cr\endpmatrix
+\pmatrix \kappa M_2\hfill& \kappa M_3\hfill\cr
\kappa M_4\hfill& \kappa^2M_5\hfill\cr\endpmatrix,
$$
where, for all $1\leq i\leq 5$, the linear operator, $M_i$, satisfies
$$
\|M_i\| \leq B,
$$
for some constant, $B>0$, which only depends on $R$. The result now follows by Lemma \procref{UniformNormOfK} and the fact that the norm of $Q_\kappa$ with respect to the $\kappa$-dependent H\"older norms, $\|\cdot\|_{k,\alpha,\kappa,\opin}$ and $\|\cdot\|_{k+1,\alpha,\kappa,\opin}$, is independent of $\kappa$.\qed
\medskip
Consider now a Banach space, $E$. Given a formal power series, $\sum_{i=0}^\infty\kappa^i f_{i,\kappa}$, taking values in $E$, and a function, $f:]0,\epsilon[\rightarrow E$, the expression
$$
f\sim \sum_{i=0}^\infty\kappa^i f_{i,\kappa}
$$
will mean that for all $m\geq 0$, there exists a constant, $B_m$, such that, for all $\kappa\in]0,1]$,
$$
\left\|f - \sum_{i=0}^m\kappa^i f_{i,\kappa}\right\| \leq B_m\kappa^{m+1}.
$$
The series, $\sum_{i=0}^\infty\kappa^i f_{i,\kappa}$, will be called the {\bf norm-asymptotic} series of the function, $f$, whenever this holds.
\proclaim{Theorem \nextprocno}
\noindent For sufficiently small $\epsilon>0$, there exist functions $X:]0,\epsilon[\rightarrow C^{k+1,\alpha}(\Bbb{R},\Bbb{R}^{d+1})$ and $\phi:]0,\epsilon[\rightarrow\hat{C}^{k+1,\alpha}_\opin(\Bbb{R}\times S^d)$ such that
$$\eqalign{
X &\sim \sum_{i=0}^\infty \kappa^i X_{i,\kappa},\cr
\phi &\sim \sum_{i=0}^\infty \kappa^i \phi_{i,\kappa}.\cr}
$$
and
$$
\Psi(\kappa,X_\kappa,\phi_\kappa) = 0.
$$
Furthermore, $X$ and $\phi$ are unique in the sense that, if $X':]0,\epsilon[\rightarrow C^{k+1,\alpha}(\Bbb{R},\Bbb{R}^{d+1})$ and $\phi':]0,1[\rightarrow C^{k+1,\alpha}_\opin(\Bbb{R}\times S^d)$ are other functions with the same properties then, for sufficiently small $\kappa$, $X(\kappa)=X'(\kappa)$, and $\phi(\kappa)=\phi'(\kappa)$.
\endproclaim
\proclabel{SingularPerturbations}
\remark In particular, for sufficiently small $\kappa$, the eternal forced mean curvature flows constructed in Theorem \procref{SingularPerturbations} are all admissable, of bounded type, and even non-degenerate, in the sense that their linearised mean curvature flow operators are surjective.
\medskip
\proof Let $m$ be a positive integer and define
$$\eqalign{
\tilde{X}_{m,\kappa} &:= \sum_{i=0}^m\kappa^i X_{i,\kappa},\ \text{and}\cr
\tilde{\phi}_{m,\kappa} &:= \sum_{i=0}^m\kappa^i \phi_{i,\kappa}.\cr}
$$
Consider the asymptotic formula, \eqnref{AsymptoticSeriesOfPsi}, for $\Psi$, substituting $X_i=X_{i,\kappa}$ and $\phi_i=\phi_{i,\kappa}$ for $0\leq i\leq m$, and $X_i=0$ and $\phi_i=0$ for $i>m$. By definition of the formal series $\sum_{i=0}^\infty\kappa^i X_{i,\kappa}$ and $\sum_{i=0}^\infty\kappa^i\phi_{i,\kappa}$, there exist smooth functions $R_{1,m}$ and $R_{2,m}$ such that
$$\eqalign{
\Pi\hat{\Psi}_\kappa(\tilde{X}_{m,\kappa},\tilde{\phi}_{m,\kappa})
&=\kappa^{m+3}\Pi\big(R_{1,m}(\gamma,D_t\gamma,x,J^1X_{0,\kappa}(t),...,J^1X_{m,\kappa}(t),\cr
&\qquad\qquad \kappa^4D_t\phi_{0,\kappa}(t,x),...,\kappa^4D_t\phi_{m,\kappa}(t,x),\cr
&\qquad\qquad\qquad\qquad J_x^2\phi_{0,\kappa}(t,x),...,J_x^2\phi_{m,\kappa}(t,x))\big),\ \text{and}\cr
\Pi^\perp\hat{\Psi}_\kappa(\tilde{X}_{m,\kappa},\tilde{\phi}_{m,\kappa})
&=\kappa^{m+1}\Pi^\perp\big(R_{2,m}(\gamma,D_t\gamma,x,J^1X_{0,\kappa}(t),...,J^1X_{m,\kappa}(t),\cr
&\qquad\qquad \kappa^4D_t\phi_{0,\kappa}(t,x),...,\kappa^4D_t\phi_{m,\kappa}(t,x),\cr
&\qquad\qquad\qquad\qquad J_x^2\phi_{0,\kappa}(t,x),...,J_x^2\phi_{m,\kappa}(t,x))\big).\cr}
$$
Recall again that composition by a given smooth function defines smooth maps between H\"older spaces which send bounded sets to bounded sets. Thus, since $R_{1,m}$ and $R_{2,m}$ are themselves smooth functions defined over finite-dimensional manifolds, it follows by the uniform bounds obtained in Lemma \procref{UniformBoundsOnCoefficients} that there exists a positive constant, $B_m$, such that, for all $\kappa$,
$$
\|T_\kappa \hat{\Psi}_\kappa(\tilde{X}_{m,\kappa},\tilde{\phi}_{m,\kappa})\|_{0,\alpha,\opin,\kappa}\leq B_m\kappa^{m+1},
$$
and the result now follows by the implicit function theorem (c.f. \cite{Rudin}).\qed
\newhead{Concentration}[Concentration]
\newsubhead{Another $\kappa$-dependent rescaled time variable}[AnotherRescaledTimeVariable]
It remains to prove the converse of Theorem \procref{SingularPerturbations}. Thus, let $F_\kappa$, $\Cal{F}_\kappa$ and $\Cal{E}_\kappa(T^{d+1})$ be defined as in Section \subheadref{EternalForcedMeanCurvatureFlows}, and, for all small $\kappa$, let $[e_{\kappa,t}]:\Bbb{R}\rightarrow\Cal{E}_\kappa(T^{d+1})$ be an admissable, eternal forced mean curvature flow of bounded type with forcing term $F_\kappa$. Given $\kappa>0$, a smooth curve, $\gamma:]a,b[\rightarrow\Bbb{R}^{d+1}$, and a smooth function, $\psi\in\hat{C}^\infty(]a,b[\times S^d)$, consider the family of embeddings, $\tilde{\Phi}(\kappa,\gamma,\psi):]a,b[\times S^d\rightarrow\Bbb{R}^{d+1}$, defined by
$$
\tilde{\Phi}(\kappa,\gamma,\psi)(r,x) := \gamma(r) + \kappa(1 + \psi(r,x))x.\eqnum{\nexteqnno[DefinitionOfTildePhi]}
$$
It turns out that the time variable denoted here by $r$ actually lies on a scale intermediate between those of $t$ and $s$. We therefore denote $t:=\kappa^2r$, and $s:=\kappa^{-2}r$, so that $D_t=\kappa^{-2}D_r$ and $D_s=\kappa^2 D_r$. However, we continue to use the $\kappa$-dependent H\"older norms, $\|\cdot\|_{k,\alpha,\kappa}$ and $\|\cdot\|_{k,\alpha,\kappa,\opin}$, which were defined in Section \subheadref{TheKappaDependentRescaledTimeVariable} using the operators $D_s$ and $\delta_s^\alpha$.

\proclaim{Lemma \nextprocno}
\noindent There exists $\kappa_0>0$ with the property that for all $\kappa<\kappa_0$, there exists a unique smooth curve, $\gamma_\kappa:\Bbb{R}\rightarrow\Bbb{T}^{d+1}$, and a unique smooth function, $\psi_\kappa\in\hat{C}^\infty(\Bbb{R}\times S^d)$, such that
$$
[e_{\kappa,r}] = [\tilde{\Phi}(\kappa,\gamma_\kappa,\psi_\kappa)].
$$
Furthermore, for all $(k,\alpha)$,
$$\eqalign{
\kappa^{-1}\|D_s\gamma_{\kappa}\|_{k,\alpha,\kappa}&\rightarrow 0,\ \text{and}\cr
\|\psi_{\kappa}\|_{k,\alpha,\kappa,\opin}&\rightarrow 0,\cr}
$$
as $\kappa$ tends to $0$.
\endproclaim
\proof We identify the flow, $[e_{\kappa,r}]$, with its lifting in $\Bbb{R}^{d+1}$. For $\kappa>0$ and $r_\kappa\in\Bbb{R}$, define the flow, $[\tilde{e}_{\kappa,s}]$, by
$$
\tilde{e}_{\kappa,s}(x) := \kappa(e_{\kappa,\kappa^2 s + r_\kappa}(x) - x_\kappa),
$$
where $x_\kappa$ is some interior point of $[e_{\kappa,r_\kappa}]$. Observe that this is an admissable, eternal forced mean curvature flow of bounded type with forcing term
$$
1 - \kappa^2 f(x_\kappa + \kappa^{-1} x).
$$
By the compactness result of Theorem \procref{Compactness}, there exists $[\tilde{e}_{0,s}]$ towards which $[\tilde{e}_{\kappa,s}]$ subconverges in the $C^k_\oploc$-sense for all $k$ as $\kappa$ tends to $0$. Furthermore, $[\tilde{e}_{0,s}]$ is an admissable, eternal forced mean curvature flow of bounded type with constant forcing term $1$, so that, by the Hopf Theorem for flows (Theorem \procref{Hopf}), it coincides up to translation with the constant flow
$$
\tilde{e}_{0,s}(x) = x.
$$
Since this limit is unique, it follows that, up to translation, the entire family, $([\tilde{e}_{\kappa,s}])_{\kappa>0}$, converges in the $C^k_\oploc$-sense for all $k$ to $[\tilde{e}_{0,s}]$. By the implicit function theorem, there therefore exists, for all sufficiently small $\kappa$, a unique curve, $\tilde{\gamma}_{\kappa}:]-1,1[\rightarrow\Bbb{R}^{d+1}$, and a unique smooth function, $\tilde{\psi}_{\kappa}\in\hat{C}^\infty(]-1,1[\times S^d)$, such that
$$
[\tilde{e}_{\kappa,s}]|_{]-1,1[} = [\tilde{\gamma}_\kappa(s) + (1 + \tilde{\psi}_\kappa(s,x))x].
$$
Furthermore, for all $(k,\alpha)$,
$$\eqalign{
\|D_s\tilde{\gamma}_\kappa\|_{k,\alpha,\kappa} &\rightarrow 0\ \text{and}\cr
\|\tilde{\psi}_\kappa\|_{k,\alpha,\kappa,\opin} &\rightarrow 0,\cr}
$$
as $\kappa$ tends to $0$. Rescaling and projecting down to the torus yields a curve $\gamma_\kappa:]r_\kappa-\kappa^2,r_\kappa+\kappa^2[\rightarrow T^{d+1}$ and a smooth function, $\psi_\kappa\in\hat{C}^\infty(]r_\kappa-\kappa^2,r_\kappa+\kappa^2[\times S^d)$, such that
$$
[e_\kappa]|_{]r_\kappa-\kappa^2,r_\kappa + \kappa^2[} = [\gamma_\kappa(r) + \kappa(1 + \psi_\kappa(r,x))x].
$$
Furthermore, for all $(k,\alpha)$,
$$\eqalign{
\kappa^{-1}\|D_s\gamma_\kappa\|_{k,\alpha,\kappa} &\rightarrow 0,\ \text{and}\cr
\|\psi_\kappa\|_{k,\alpha,\kappa,\opin} &\rightarrow 0,\cr}
$$
as $\kappa$ tends to $0$. Finally, since the family $(r_\kappa)_{\kappa>0}$ is arbitrary, this convergence is uniform over $\Bbb{R}\times S^d$, and the result follows.\qed
\newsubhead{Bootstrapping}[Bootstrapping]
Consider the function, $\tilde{\Psi}:]0,\infty[\times J\rightarrow\Bbb{R}$, defined such that, for all $\kappa\in]0,\infty[$, for all smooth functions, $\gamma:\Bbb{R}\rightarrow T^{d+1}$ and $\psi:\Bbb{R}\times S^d\rightarrow\Bbb{R}$, and for any point $(r,x)\in\Bbb{R}\times S^d$,
$$\eqalign{
\tilde{\Psi}(\kappa,J^1_t\gamma(r),J^2_\opin\psi(r,x))
&:=\langle D_r\tilde{\Phi}(\kappa,\gamma,\phi),N(J_x^1\psi(r,x))\rangle\cr
&\qquad\qquad + \frac{1}{\kappa}H(J_x^2\psi(r,x)) - F_\kappa(\tilde{\Phi}(\kappa,\gamma,\psi(r,x)).\cr}\eqnum{\nexteqnno[DefinitionOfTildePsi]}
$$
As before, $\tilde{\Psi}$ is a smooth function defined over the finite-dimensional manifold, $]0,\infty[\times J^1$. Likewise, by \eqnref{ForcedMCFEquation}, $\tilde{\Psi}$ has been constructed precisely so that, given $\kappa>0$, and smooth functions, $\gamma:\Bbb{R}\rightarrow T^{d+1}$ and $\phi:\Bbb{R}\times S^d\rightarrow\Bbb{R}$, the family, $\tilde{\Phi}(\kappa,X,\psi)$, is an eternal forced mean curvature flow with forcing term $F_\kappa$ if and only if $\tilde{\Psi}(\kappa,J^1\gamma,J^1_\opin\psi)$ vanishes at every point, $(r,x)$, of $\Bbb{R}\times S^d$.
\proclaim{Lemma \nextprocno}
\noindent For all $(k,\alpha)$, there exists $B_{k,\alpha}>0$ such that, for all $\kappa\in]0,1]$,
$$\eqalign{
\|D_r\gamma_\kappa + \kappa^2\nabla f(\gamma_\kappa(r))\|_{k,\alpha,\kappa} &\leq B_{k,\alpha}\kappa^3,\ \text{and}\cr
\|\psi_\kappa\|_{k,\alpha,\kappa,\opin} &\leq B_{k,\alpha}\kappa^2.\cr}
$$
\endproclaim
\proclabel{Bootstrapping}
\proof By \eqnref{TaylorSeriesOfH}, \eqnref{TaylorSeriesOfN} and \eqnref{DefinitionOfTildePsi},
$$\eqalign{
\tilde{\Psi} &= \frac{1}{\kappa^2}\langle D_s\gamma_\kappa,x\rangle - \frac{1}{\kappa^2}\langle D_s\gamma_\kappa,\nabla\psi_\kappa\rangle + \frac{1}{\kappa^2}\langle D_s\gamma-\kappa,R_1(x,J_x^1\psi_\kappa)\rangle\cr
&\qquad + \frac{1}{\kappa}D_s\psi_\kappa - \frac{1}{d\kappa}(d + \Delta)\psi_\kappa + \frac{1}{\kappa}R_2(x,D_s\psi_\kappa,J_x^2\psi_\kappa)\cr
&\qquad\qquad + \kappa f(\gamma_\kappa) + \kappa^2\langle\nabla f(\gamma_\kappa),x\rangle + \kappa^2\psi_\kappa\langle\nabla f(\gamma_\kappa),x\rangle + \kappa^3 R_3(\kappa,x,\gamma_\kappa,\psi_\kappa),\cr}
$$
where $R_1$, $R_2$ and $R_3$ are smooth functions of their arguments and $R_1$, $R_2$ are of order at least $2$ in $(D_s\psi_\kappa,J_x^2\psi_\kappa)$. Projecting onto $\Cal{H}_1^\perp$ yields
$$\eqalign{
D_s\psi_\kappa + \langle\kappa^{-1}D_s\gamma_\kappa,\nabla\psi_\kappa\rangle + \frac{1}{d}(d+\Delta)\psi_\kappa
&=\Pi^\perp R_4(x,\kappa^{-1}D_s\gamma_\kappa,D_s\psi_\kappa,J_s^2\psi_\kappa)\cr
&\qquad\qquad + \kappa^2\Pi^\perp R_5(\kappa,x,\gamma_\kappa,\psi_\kappa),\cr}
$$
where $R_4$ and $R_5$ are smooth functions of their arguments and $R_4$ is of order at least $2$ in $(D_s\psi_\kappa,J_x^2\psi_\kappa)$. Since the operator $D_s+\frac{1}{d}(d+\Delta)$ maps $\hat{C}^{k+1,\alpha}(\Bbb{R}\times S^d)$ invertibly onto $\hat{C}^{k,\alpha}(\Bbb{R}\times S^d)$, it follows from the estimates already obtained for $D_s\gamma_\kappa$ that the operator on the left hand side is uniformly bounded below for sufficiently small $\kappa$. For all $(k,\alpha)$, there therefore exists $B>0$ such that
$$
\|\psi_{\kappa}\|_{k,\alpha,\kappa,\opin} \leq B\|\psi_\kappa\|_{k,\alpha,\kappa,\opin}^2 + B\kappa^2,
$$
and since $\|\psi_\kappa\|_{k,\alpha,\kappa,\opin}$ tends to $0$ as $\kappa$ tends to $0$, the second assertion now follows. Projecting onto $\Cal{H}_1$ now yields
$$\eqalign{
(\opId + M(J_x^1\psi_\kappa))D_t\gamma_\kappa + \kappa^2\nabla f(\gamma_\kappa)
&=\frac{1}{\kappa}\Pi R_6(x,D_s\psi_\kappa,J_x^2\psi_\kappa) + \kappa^2\Pi R_7(x,J_x^0\psi_\kappa)\cr
&\qquad\qquad +\kappa^3\Pi R_8(\kappa,x,\gamma_\kappa,\psi_\kappa),}
$$
where $M$, $R_6$, $R_7$ and $R_8$ are smooth functions of their arguments, $R_6$ has order at least $2$ in $(D_s\psi_\kappa,J_x^2\psi_\kappa)$ and $M$ and $R_7$ have order $1$ in $J_x\phi$, and the first assertion now follows by the bounds obtained on $\psi_\kappa$.\qed\goodbreak
\newsubhead{Recovering the Flow}[RecoveringTheFlow]
It follows from Theorem \procref{ConcentrationEllipticCase} that, as $\kappa$ tends to $0$, the two end-points of the flow, $[e_{\kappa,t}]$, collapse onto critical points of $f$. In fact (c.f. \cite{SmiEC} and \cite{Ye}), for sufficiently small $\kappa$, they coincide with those Alexandrov-embeddings constructed via the elliptic analogue of Theorem \procref{SingularPerturbations}. We now show that the entire flow in fact collapses onto a complete gradient flow of $f$ as the parameter, $\kappa$, tends to $0$.
\proclaim{Lemma \nextprocno}
\noindent Let $M$ be a compact manifold. If $f:M\rightarrow\Bbb{R}$ is of Morse-Smale type, then there exists $C>0$ with the following property. If $\epsilon>0$ and if $\gamma:\Bbb{R}\rightarrow M$ is a $C^1$ curve such that
$$\eqalign{
\msup_{t}\|\dot{\gamma}(t) + \nabla f(\gamma(t))\| &< \epsilon,\cr
\mlimsup_{t\rightarrow-\infty} d(\gamma(t),p_-) &< \epsilon,\ \text{and}\cr
\mlimsup_{t\rightarrow+\infty} d(\gamma(t),p_+) &< \epsilon,\cr}
$$
where $p_-$ and $p_+$ are critical points of $f$ and if
$$
\opIndex(p_-) - \opIndex(p_+) = 1,
$$
then there exists a complete integral curve $\gamma_0:\Bbb{R}\rightarrow M$ of $-\nabla f$ such that
$$\eqalign{
\mlim_{t\rightarrow-\infty}\gamma_0(t) &= p_-,\cr
\mlim_{t\rightarrow+\infty}\gamma_0(t) &= p_+,\cr
\msup_{t}\|\gamma(t) - \gamma_0(t)\| &< C\epsilon,\ \text{and}\cr
\int_{-\infty}^{+\infty}\langle\dot{\gamma}_0,\gamma(t)-\gamma_0(t)\rangle dt &=0.\cr}
$$
\endproclaim
\proclabel{RecoveringTheFlow}
\remark Since $F$ is of Morse-Smale type, $\dot{\gamma}_0$, decays exponentially at $\pm\infty$, so that the above integral is always well-defined.
\medskip
\proof We first show that for all $\delta>0$ there exists $\epsilon_0$ such that for $\epsilon<\epsilon_0$, there exists a complete integral curve, $\gamma_0:\Bbb{R}\rightarrow M$, of $-\nabla f$ such that
$$\eqalign{
\mlim_{t\rightarrow-\infty}\gamma_0(t) &= p_-,\cr
\mlim_{t\rightarrow+\infty}\gamma_0(t) &= p_+,\ \text{and}\cr
\msup_{t}\|\gamma(t) - \gamma_0(t)\| &< \delta.\cr}
$$
Indeed, choose $\delta>0$. By non-degeneracy and compactness, we may suppose that if $p$ is a critical point of $f$ and if $\eta:[0,\infty[\rightarrow M$ is a gradient flow of $f$ starting at some point of $\partial B_{2\delta}(p)$, then
$$
\eta([0,\infty[) \minter B_{3\delta}(q) =\emptyset
$$
for every other critical point, $q$, of $f$ such that $\opIndex(q)\geq\opIndex(p)$. For every critical point, $p$, of $f$, let $S_{2\delta}(p)$ be the intersection of the stable manifold of $p$ with $\partial B_{2\delta}(p)$. By compactness again, for all $\delta'>0$, there exists $T>0$ such that if $\eta:[0,\infty[\rightarrow M$ is a gradient flow of $f$ starting at some point of $\partial B_{2\delta}(p)$ such that
$$
d(\eta(0),S_{2\delta}(p)) > \delta',
$$
then
$$
\eta([0,T])\minter B_\delta(q) \neq \emptyset,
$$
for some other critical point $q$ of $f$, which, in particular, has Morse index strictly less than that of $p$. However, for sufficiently small $\epsilon$, if $\gamma(t_0)\in\partial B_{2\delta}(p)$ is such that $\dot{\gamma}(t_0)$ points outward from $B_{2\delta}(p)$ or is tangent to $\partial B_{2\delta}(p)$, then, on the one hand, $d(\gamma(t_0),S_{2\delta}(p))>\delta'$, and, on the other,
$$
d(\eta(t),\gamma(t)) < \delta\ \forall\ t_0\leq t\leq t+T,
$$
where $\eta$ is the gradient flow of $f$ such that $\eta(t_0)=\gamma(t_0)$. Consequently, if $\gamma$ leaves $B_{2\delta}(p)$, then, after a time at most $T$, it must enter $B_{2\delta}(q)$ for some critical point, $q$, of $f$ of Morse index strictly less than that of $p$, after which it cannot return to $B_{2\delta}(p)$.
\par
For $\epsilon$ sufficiently small, $\gamma$ starts inside $B_{2\delta}(p_-)$. By the above discussion, within a time at most $T$ after leaving $B_{2\delta}(p_-)$, the curve, $\gamma$, enters $B_{2\delta}(p_+)$ from which it does not leave. Furthermore, there exists a gradient flow $\eta$ of $f$ from some point of $B_{2\delta}(p_-)$ to some point of $B_{2\delta}(p_+)$ from which $\gamma$ remains at a distance of at most $\delta$ as it travels between these two balls. This proves the assertion.
\par
Now define $X:=\gamma(t)-\gamma_0(t)$. In particular,
$$
\|X\|_{C^0} < \delta.
$$
Since $f$ is of Morse-Smale type, the operator $L:=\partial_t-\opHess(f)(\gamma(t))$ is Fredholm and surjective. Furthermore, since $\opIndex(p_-)-\opIndex(p_+)=1$, by the Atiyah-Patodi-Singer index theorem, its kernel is $1$-dimensional and is spanned by $\dot{\gamma}_0$. Observe that, upon replacing $\gamma_0(t)$ with $\gamma_0(s+t)$ for a suitable value of $s$, we may suppose that
$$
\int_{-\infty}^{+\infty}\langle\dot{\gamma}_0(t),X(t)\rangle dt = 0,
$$
so that $X$ is an element of the orthogonal complement of the kernel of $P$. However,
$$\eqalign{
\|LX\|_0 &= \|\dot{X}(t) - \opHess(f)(\gamma_0(t))X(t)\|_0\cr
&= \|\dot{\gamma}(t) - \dot{\gamma}_0(t) - \opHess(f)(\gamma_0(t))X(t)\|_0\cr
&\leq \epsilon + \|\nabla f(\gamma(t)) - \nabla f(\gamma_0(t)) - \opHess(f)(\gamma_0(t))X(t)\|_0\cr
&\leq \epsilon + B_1\|X(t)\|^2_0,\cr}
$$
for some constant $B_1$, and since the restriction of $L$ to the orthogonal complement of $\dot{\gamma}_0$ is invertible,
$$
\|\dot{X} - \opHess(f)(\gamma(t))X\|_0 \geq \frac{1}{B_2}\|X\|_0,
$$
for a suitable constant $B_2$. Combining these relations yields
$$
\|X\|_0\left(1 - B_1B_2\|X\|_0\right) \leq B_2\epsilon,
$$
and since $\|X\|_0<\delta$, the result follows.\qed
\newsubhead{Asymptotic series}[AsymptoticSeriesAgain]
Finally, recall the function, $\Psi$, defined in Section \subheadref{AsymptoticSeries}. It follows from Lemmas \procref{Bootstrapping} and \procref{RecoveringTheFlow} that for all sufficiently small $\kappa$, there exists $X_\kappa\in C^\infty_\opbdd(\Bbb{R},\Bbb{R}^{d+1})$ and $\phi_\kappa\in\hat{C}^\infty_\opbdd(\Bbb{R}\times S^d)$ such that
$$
[e_{\kappa,t}] = \Psi(\kappa,X_\kappa,\phi_\kappa).
$$
Furthermore, for all $(k,\alpha)$, there exists $B_{k,\alpha}>0$ such that, for all $\kappa\in]0,1]$,
$$\eqalign{
\|X_\kappa\|_{k,\alpha,\kappa}&\leq B_{k,\alpha},\cr
\|D_tX_\kappa\|_{k,\alpha,\kappa}&\leq B_{k,\alpha},\ \text{and}\cr
\|\phi_\kappa\|_{k,\alpha,\kappa,\opin}&\leq B_{k,\alpha}.\cr}
$$
Consider now the formal solutions, $\sum_{k=0}^\infty\kappa^k X_{m,\kappa}$ and $\sum_{k=0}^\infty\kappa^k\phi_{m,\kappa}$, constructed in Lemma \procref{FormalSolutions}.
\proclaim{Lemma \nextprocno}
\noindent The family $(X_\kappa,\phi_\kappa)$ satisfies
$$\eqalign{
X_\kappa &\sim \sum_{k=0}^\infty \kappa^k X_{k,\kappa},\ \text{and}\cr
\phi_\kappa &\sim \sum_{k=0}^\infty \kappa^k \phi_{k,\kappa}.}
$$
\endproclaim
\proclabel{FlowsAreAsymptoticToFormalSolution}
\proof This is proven by induction. Indeed, denote
$$\eqalign{
\tilde{X}_\kappa &:= \kappa^{-(m+1)}\left(X_\kappa - \sum_{k=0}^m\kappa^k X_{k,\kappa}\right),\ \text{and}\cr
\tilde{\phi}_\kappa &:= \kappa^{-(m+1)}\left(\phi_\kappa -\sum_{k=0}^m\kappa^k \phi_{k,\kappa}\right).\cr}
$$
Suppose, by induction, that, for all $(k,\alpha)$, there exists $B_{k,\alpha}$ such that, for all $\kappa\in]0,1]$,
$$\eqalign{
\|\tilde{X}_\kappa\|_{k,\alpha,\kappa}&\leq B_{k,\alpha},\cr
\|D_t\tilde{X}_\kappa\|_{k,\alpha,\kappa}&\leq B_{k,\alpha},\ \text{and}\cr
\|\tilde{\phi}_\kappa\|_{k,\alpha,\kappa}&\leq B_{k,\alpha}.\cr}
$$
Consider the asymptotic formula, \eqnref{AsymptoticSeriesOfPsi}, for $\Psi$, substituting $X_k=X_{k,\kappa}$, for $0\leq k\leq m$, $X_{m+1}=\tilde{X}_\kappa$, and $X_k=0$, for $k\geq m+2$; and substituting $\phi_k=\phi_{k,\kappa}$, for $0\leq k\leq m$, $\phi_{m+1}=\tilde{\phi}_\kappa$ and $\phi_k=0$, for $k\geq m+2$. Projecting the $(m+2)$'nd order term orthogonally onto $\hat{\Cal{H}}$, we see that there exists a smooth function, $R_1$, such that
$$\eqalign{
P_\kappa\tilde{\phi}_\kappa
&=-\langle\nabla f(\gamma(t)),X_{m,\kappa}\rangle - \langle\opHess(f)(\gamma(t))X_{m-1,\kappa},X_{0,\kappa}\rangle\cr
&\qquad\qquad -\Pi^\perp\big(\Psi_{m+1}(\gamma,D_t\gamma,x,J_t^1 X_{0,\kappa},...,J_t^1X_{m-2,\kappa},\cr
&\qquad\qquad\qquad\qquad \kappa^4D_t\phi_{0,\kappa},...,\kappa^4D_t\phi_{m-1,\kappa},J_x^2\phi_{0,\kappa},...,J_x^2\phi_{m-1,\kappa})\big)\cr
&\qquad\qquad -\kappa \Pi^\perp\big(R_1(\gamma,D_t\gamma,x,J_t^1 X_{0,\kappa},...,J_t^1X_{m,\kappa},J_t^1\tilde{X}_\kappa\cr
&\qquad\qquad\qquad\qquad \kappa^4D_t\phi_{0,\kappa},...,\kappa^4D_t\phi_{m,\kappa},\kappa^4D_t\tilde{\phi}_\kappa,
J_x^2\phi_{0,\kappa},...,J_x^2\phi_{m,\kappa},J_x^2\tilde{\phi}_\kappa)\big).\cr}
$$
However, for all $(k,\alpha)$, the operator $P_\kappa$ defines an invertible linear map from $\hat{C}^{k+1,\alpha}_{\opin}(\Bbb{R}\times S^d)$ into $\hat{C}^{k,\alpha}_{\opin}(\Bbb{R}\times S^d)$. It therefore follows from the definition, \eqnref{InductiveDefinitionOfPhi}, of $\phi_{k+1,\kappa}$ that there exists a constant, $B>0$, such that, for all $\kappa\in]0,1]$,
$$
\|\tilde{\phi}_\kappa - \phi_{k+1,\kappa}\|_{k,\alpha,\kappa} \leq B\kappa.
$$
Denote now
$$
\tilde{\phi}_\kappa' := \kappa^{-1}\left(\tilde{\phi}_\kappa - \phi_{k+1,\kappa}\right).
$$
Consider again the asymptotic formula, \eqnref{AsymptoticSeriesOfPsi}, for $\Psi$, substituting again $X_k=X_{k,\kappa}$, for $0\leq k\leq m$, $X_{m+1}=\tilde{X}_\kappa$ and $X_k=0$, for $k\geq m+2$; but substituting this time $\phi_k=\phi_{k,\kappa}$, for $0\leq k\leq m+1$, $\phi_{m+2}=\tilde{\phi}_\kappa'$ and $\phi_k=0$, for $k\geq m+3$. Projecting the $(m+4)$'th order term orthogonally onto $\Cal{H}_1$, we see that there exists a smooth function, $R_2$, such that
$$\eqalign{
L\tilde{X}
&=-\Pi\big(\Psi_{m+1}(\gamma,D_t\gamma,x,J_t^1 X_{0,\kappa},...,J_t^1X_{m,\kappa},\cr
&\qquad\qquad \kappa^4D_t\phi_{0,\kappa},...,\kappa^4D_t\phi_{m+1,\kappa},J_x^2\phi_{0,\kappa},...,J_x^2\phi_{m+1,\kappa})\big)\cr
&\qquad -\kappa \Pi\big(R_2(\gamma,D_t\gamma,x,J_t^1 X_{0,\kappa},...,J_t^1X_{m,\kappa},J_t^1\tilde{X}_\kappa\cr
&\qquad\qquad \kappa^4D_t\phi_{0,\kappa},...,\kappa^4D_t\phi_{m+1,\kappa},\kappa^4 D_t\tilde{\phi}'_\kappa,J_x^2\phi_{0,\kappa},...,J_x^2\phi_{m+1,\kappa},J_x^2\tilde{\phi}'_\kappa)\big),\cr}
$$
However, by Lemma \procref{UniformNormOfK}, for all $(k,\alpha)$, the norm of $K$ is uniformly bounded independent of $\kappa\in]0,1]$. It therefore follows by the definition, \eqnref{InductiveDefinitionOfX}, of $X_{k+1,\kappa}$ that there exists a constant, $B>0$, such that
$$
\|\tilde{X} - X_{k+1,\kappa}\|_{k,\alpha} \leq B\kappa.
$$
The result now follows by induction.\qed
\inappendicestrue
\global\headno=0
\medskip
\newhead{The weakly smooth category}[TheWeaklySmoothCategory]
Let $M$ be a $(d+1)$-dimensional riemannian manifold, and let $\Cal{E}(M)$ be defined as in the introduction. Charts of $\Cal{E}(M)$ are constructed as follows. For $[e]\in\Cal{E}(M)$, let $N_e:S^d\rightarrow TM$ denote the unit, normal vector field over $e$ which is compatible with the orientation, and consider the map, $E_e:S^d\times\Bbb{R}\rightarrow M$, given by
$$
E_e(x,t) := \opExp_{e(x)}(tN_e(x)),
$$
where $\opExp$ is the exponential map of $M$. Let $\epsilon_e>0$ be such that the restriction of $E_e$ to $S^d\times]-\epsilon_e,\epsilon_e[$ is an immersion (c.f. Proposition $3.2$ of \cite{MaximoNunesSmith}), and define the open subset, $\Cal{U}_e$, of $C^\infty(S^d)$ and the mapping $\hat{\Phi}_e:\Cal{E}_e\rightarrow C^\infty(S^d,M)$ by
$$\eqalign{
\Cal{U}_e &:= \left\{\phi\in C^\infty(\Sigma)\ |\ \|\phi\|_{C^0}<\epsilon_e\right\},\ \text{and}\cr
\hat{\Phi}_e(\phi)(x) &:= E_e(x,\phi(x)).\cr}
$$
Upon reducing $\epsilon_e$ if necessary, we may suppose that $\hat{\Phi}_e(\phi)$ is an Alexandrov-embedding for all $\phi\in\Cal{U}_e$ so that $\hat{\Phi}_e$ projects down to a mapping, $\Phi_e$, which is in fact a homeomorphism from $\Cal{U}_e$ into an open subset, $\Cal{V}_e$, of $\Cal{E}(M)$ (c.f. Propositions $3.3$ and $3.4$ of \cite{MaximoNunesSmith}). The triplet $(\Phi_e,\Cal{U}_e,\Cal{V}_e)$ is called the {\bf graph chart} of $\Cal{E}(M)$ about $e$. The set of all graph charts consitutes an atlas of $\Cal{E}(M)$ all of whose transition maps are homeomorphisms, thereby furnishing $\Cal{E}(M)$ with the structure of a topological manifold modelled on $C^\infty(S^d)$.
\par
This atlas in fact furnishes $\Cal{E}(M)$ with the structure of a smooth-tame Frechet manifold (c.f. \cite{HamiltonII} for a splendid introduction to this theory). However, since this framework is highly technical, we prefer to avoid it whenever possible. It is for this reason that we introduce the following formalism of weakly smooth manifolds which, despite its simplicity, possesses all the structure required for the development of the formal parts of the theory. Thus, consider first an open subset, $\Cal{U}$, of $C^\infty(S^d)$. Given a finite-dimensional manifold, $X$, and a function, $c:X\rightarrow\Cal{U}$, the function, $\tilde{c}:X\times S^d\rightarrow M$, is defined by $\tilde{c}(x,y):=c(x)(y)$. The function, $c$, is then said to be {\bf strongly smooth} whenever the function, $\tilde{c}$, is smooth. Now consider another open subset, $\Cal{U}'$, of $C^\infty(S^d)$. A function, $\Phi:\Cal{U}\rightarrow\Cal{U}'$, is said to be {\bf weakly smooth} whenever composition by $\Phi$ sends strongly smooth functions continuously into strongly smooth functions. Observe, in particular, that every weakly smooth function is also continuous. It is now a straightforward matter to show that the transition maps between graph charts of $\Cal{E}(M)$ are weakly smooth, so that $\Cal{E}(M)$ carries the structure of a weakly smooth manifold. In particular, a function $c:X\rightarrow\Cal{E}(M)$ is now said to be {\bf strongly smooth} whenever it is strongly smooth in every graph chart.
\par
The formalism of weakly smooth manifolds possesses considerable structure. For example, tangent vectors are well-defined, as is the tangent bundle, which is also a weakly smooth manifold, and so on. The main important result not admitted by this formalism is the inverse function theorem. However, this presents no problem in the present context, since ellipticity and hypo-ellipticity allow us to work locally in every graph chart as if it were a H\"older space, where the inverse function theorem for Banach manifolds can then be applied. For a more thorough discussion, we refer the reader to \cite{RosSmiDT}.
\newhead{Bibliography}[Bibliography]
{\leftskip = 5ex \parindent = -5ex
\leavevmode\hbox to 4ex{\hfil \cite{AgmonNirenberg}}\hskip 1ex{Agmon S., Nirenberg L., Lower bounds and uniqueness theorems for solutions of differential equations in a Hilbert space, {\sl Comm. Pure Appl. Math.}, {\bf 20}, (1967), 207–-229}
\medskip
\leavevmode\hbox to 4ex{\hfil \cite{Alexander}}\hskip 1ex{Alexander S., Locally convex hypersurfaces of negatively curved spaces, {\sl Proc. Amer. Math. Soc.}, {\bf 64}, (1977), no. 2, 321--325}
\medskip
\leavevmode\hbox to 4ex{\hfil \cite{AlmgrenSimon}}\hskip 1ex{Almgren F. J., Simon L., Existence of embedded solutions of Plateau's problem, {\sl Ann. Scuola Norm. Sup. Pisa Cl. Sci.}, (\bf 6), (1979), no. 3, 447--495}
\medskip
\leavevmode\hbox to 4ex{\hfil \cite{Andrews}}\hskip 1ex{Andrews B., Non-collapsing in mean-convex mean curvature flow, {\sl Geometry \& Topology}, {\bf 16}, (2012), 1413--1418}
\medskip
\leavevmode\hbox to 4ex{\hfil \cite{AndrewsLangfordMcCoy}}\hskip 1ex{Andrews B., Langford M., McCoy J., Non-collapsing in fully nonlinear curvature flows, {\sl Annales de l'Institut Henri Poincar\'e - Analyse non lin\'eaire}, {\bf 30}, (2013), no. 1, 23--32}
\medskip
\leavevmode\hbox to 4ex{\hfil \cite{GrueterJost}}\hskip 1ex{Gr\"uter M., Jost J., On embedded minimal disks in convex bodies. {\sl Ann. Inst. H. Poincar\'e, Anal. Non Lin\'eaire}, {\bf 3}, (1986), no. 5, 345--390}
\medskip
\leavevmode\hbox to 4ex{\hfil \cite{Hamilton}}\hskip 1ex{Hamilton R. S., Convex hypersurfaces with pinched second fundamental form, {\sl Comm. Anal. Geom.}, {\bf 2}, (1994), no. 1, 167--172}
\medskip
\leavevmode\hbox to 4ex{\hfil \cite{HamiltonII}}\hskip 1ex{Hamilton R. S., The inverse function theorem of Nash and Moser, {\sl Bull. Amer. Math. Soc.}, {\bf 7}, (1982), no. 1, 65--222}
\medskip
\leavevmode\hbox to 4ex{\hfil \cite{Mantegazza}}\hskip 1ex{Mantegazza C., {\sl Lecture notes on mean curvature flow}, Progress in Mathematics, {\bf 290}, Springer Verlag, Basel, (2011)}
\medskip
\leavevmode\hbox to 4ex{\hfil \cite{MacSmiMT}}\hskip 1ex{Macarini L., Smith G., Morse homology of forced mean curvature flows, {\sl in preparation}}
\medskip
\leavevmode\hbox to 4ex{\hfil \cite{MaximoNunesSmith}}\hskip 1ex{M\'aximo D., Nunes I., Smith G., Free boundary minimal annuli in convex three-manifolds, to appear in {\sl J. Diff. Geom}}
\medskip
\leavevmode\hbox to 4ex{\hfil \cite{PacardXu}}\hskip 1ex{Pacard F., Xu X., Constant mean curvature spheres in Riemannian manifolds,\break {\sl Manuscripta Math.}, {\bf 128}, (2009), no. 3, 275--295}
\medskip
\leavevmode\hbox to 4ex{\hfil \cite{RobbinSalamon}}\hskip 1ex{Robbin J., Salamon D., The spectral flow and the Maslov index, {\sl Bull. London Math. Soc.}, {\bf 27}, (1995), no. 1, 1--33}
\medskip
\leavevmode\hbox to 4ex{\hfil \cite{RosSmiDT}}\hskip 1ex{Rosenberg H., Smith G., Degree theory of immersed hypersurfaces, arXiv:1010.1879}
\medskip
\leavevmode\hbox to 4ex{\hfil \cite{Rudin}}\hskip 1ex{Rudin W., {\sl Principles of mathematical analysis}, International Series in Pure \& Applied Mathematics, McGraw-Hill, (1976)}
\medskip
\leavevmode\hbox to 4ex{\hfil \cite{Schneider}}\hskip 1ex{Schneider M., Closed magnetic geodesics on closed hyperbolic Riemann surfaces, {\sl Proc. London Math. Soc.}, (2012), {\bf 105}, 424--446}
\medskip
\leavevmode\hbox to 4ex{\hfil \cite{Schwarz}}\hskip 1ex{Schwarz M., {\sl Morse homology}, Progress in Mathematics, {\bf 111}, Birkh\"auser Verlag, Basel, (1993)}
\medskip
\leavevmode\hbox to 4ex{\hfil \cite{SmiEC}}\hskip 1ex{Smith G., Constant curvature hypersurfaces and the Euler characteristic,\break arXiv:1103.3235}
\medskip
\leavevmode\hbox to 4ex{\hfil \cite{SmiEMCFI}}\hskip 1ex{Smith G., Eternal forced mean curvature flows I - a compactness result, {\sl Geom. Dedicata}, {\bf 176}, no. 1, (2014), 11--29}
\medskip
\leavevmode\hbox to 4ex{\hfil \cite{SmiEMCFII}}\hskip 1ex{Smith G., Eternal Forced Mean Curvature Flows II - Existence, arXiv:1508.05688}
\medskip
\leavevmode\hbox to 4ex{\hfil \cite{Tromba}}\hskip 1ex{Tromba A, {\sl A theory of branched minimal surfaces}, Springer Monographs in Mathematics, Springer-Verlag, Berlin, Heidelberg, (2012)}
\medskip
\leavevmode\hbox to 4ex{\hfil \cite{Weber}}\hskip 1ex{Weber J., Morse homology for the heat flow, {\sl Math. Z.}, {\bf 275}, (2013), no. 1, 1--54}
\medskip
\leavevmode\hbox to 4ex{\hfil \cite{WhiteI}}\hskip 1ex{White B., Every three-sphere of positive Ricci curvature contains a minimal embedded torus, {\sl Bull. Amer. Math. Soc.}, {\bf 21}, (1989), no. 1, 71--75}
\medskip
\leavevmode\hbox to 4ex{\hfil \cite{WhiteII}}\hskip 1ex{White B., The space of minimal submanifolds for varying Riemannian metrics, {\sl Indiana Univ. Math. J.}, {\bf 40}, (1991), no. 1, 161--200}
\medskip
\leavevmode\hbox to 4ex{\hfil \cite{Ye}}\hskip 1ex{Ye R., Foliation by constant mean curvature spheres, {\sl Pacific J. Math.}, {\bf 147}, (1991), no. 2, 381--396}
\par}
%
%
%
%
\enddocument